%
%
%

\input kt3.mac

\input macro-lagrangian.mac

\magnification=\magstephalf
\baselineskip=15 true pt

\input macro-pdftex.mac

\input bbb-plain.mac


\citenumber

\advance\vsize by 1cm


\centerline{\bf A Lagrangian approach for prescribed mass solutions of}
\centerline{\bf  cubic-quintic Schr\"odinger equations and $L^2$-supercritical problems}

\bigskip

\centerline{Silvia Cingolani${}^1$, Marco Gallo${}^2$, Kazunaga Tanaka${}^3$}

\bigskip

\settabs 20 \columns

\+&\hfil${}^1$&Dipartimento di Matematica, Universit\`{a} degli Studi di Bari Aldo Moro\cr
\+&&Via E. Orabona 4, 70125, Bari, Italy\cr
\+&&{\tt ORCID: 0000-0002-3680-9106}\cr

\smallskip

\+&\hfil${}^2$&Dipartimento di Matematica e Fisica, Universit\`{a} Cattolica del Sacro Cuore\cr
\+&&Via della Garzetta 48, 25133 Brescia, Italy\cr
\+&&{\tt ORCID: 0000-0002-3141-9598}\cr

\smallskip

\+&\hfil${}^3$&Department of Mathematics, School of Science and Engineering, Waseda University\cr
\+&&3-4-1 Ohkubo, Shijuku-ku, Tokyo 169-8555, Japan\cr
\+&&{\tt ORCID: 0000-0002-1144-1536}\cr

\bigskip

{\narrower

\noindent
{\bf Abstract.}\m  
We study the existence of radially symmetric solutions of the following 
nonlinear scalar field equations in $\R^N$ $(N \geq 2)$:
    $$(*)_m \quad 
        - \Delta u + \mu u = g(u) \quad \hbox{in}\ \R^N, 
        \quad  \half  \intRN u^2\, dx = m, 
    $$
where $g(s) \in C(\R,\R)$, $m > 0$ and $\mu \in \R$ is an unknown Lagrangian multiplier.  
We take an approach using a Lagrangian formulation of $(*)_m$:
    $$  \eqalign{
        J_m(\mu,u)&=\half\intRN\abs{\nabla u}^2\,dx -\intRN G(u)\,dx 
            +\mu\left(\half\intRN u^2\, dx-m\right) \cr
        &\in C^1((0,\infty)\times H_r^1(\R^N), \R) \cr}
    $$
and we give new general existence results through the function:
    $$  \eqalign{
        b_m:\, (0,\infty) \to \R; 
        \mu \mapsto \hbox{Mountain Pass minimax value for}\  
            (u\mapsto J_m(\mu,u)).\cr}
    $$
We will show the existence of solutions of $(*)_m$ related to local minima 
and local maxima of $b_m(\mu)$.  
As applications, we study cubic-quintic type equations and $L^2$-supercritical problems.
In particular, when $N=2,3$, 
we show new existence results of normalized solutions
without assuming global Ambrosetti-Rabinowitz type conditions, which partially improve
the preceding results due to Jeanjean [\cite[J:24]] and 
Jeanjean-Lu [\cite[JL1:26], \cite[JL3:28]].

}


\bigskip


\noindent
{\bf keywords:} nonlinear Schr\"odinger equations, standing waves, normalized solutions, 
cubic-quintic model, mass-supercritical problem, Lagrangian formulation.


\bigskip

\settabs 18\columns

{\bf Contents}

\medskip

\tabalign &\ref[Section:1]. Introduction\cr
\tabalign &\ref[Section:2]. Preliminaries\cr
\tabalign &\ref[Section:3]. Behaviors of $a(\mu)$ and $c_{-}(\mu)$\cr
\tabalign &\ref[Section:4]. $(PSPC)$ condition for $\widetilde{J}_m(\lambda,u)$\cr
\tabalign &\ref[Section:5]. Application to cubic-quintic type equations and $L^2$-supercritical problems\cr
\tabalign &\ref[Section:6]. Proofs of abstract results\cr

\tabalign&\link{References}{References}\cr

\vfil
\break


\BS{\label[Section:1]. Introduction}
In this paper we study the existence of radially symmetric
positive solutions of the following nonlinear scalar field
equations in $\R^N$ $(N \geq 2)$:
    $$  (\ast)_m\ \left\{
        \eqalign{
            - &\Delta u + \mu u = g(u) \quad \hbox{in } \R^N, \cr
            &\half  \int_{\R^N} u^2\, dx  = m, \cr
            &u \in H^1_r(\R^N), \cr}
        \right.
        $$
where $g(s) \in C(\R,\R)$, $m>0$ is a given constant,
and $\mu \in \R$ is an unknown Lagrange multiplier.

Solutions of $(*)_m$ can be characterized as critical
points of the constraint problem:
    $$  \calI(u) = \half  \int_{\R^N} \abs{\nabla u}^2\, dx 
        - \int_{\R^N} G(u)\, dx : \, \calS_m \to \R,
        \eqno\label[1.1]
    $$
where $\calS_m = \{ u \in H^1_r(\R^N) \,;\, {1\over 2} \int_{\R^N} u^2\, dx = m \}$  
and $G(s) = \int_0^s g(t)\,dt$.
Many researchers study $(*)_m$ through $\calI(u)$. See  
Cazenave-Lions [\cite[CL:12]], Shibata [\cite[Sh:38]], Ruppen [\cite[Ru:35]], Stuart [\cite[St:40]], 
Hirata-Tanaka [\cite[HT:22]], Jeanjean-Lu [\cite[JL2:27]], Ikoma-Miyamoto [\cite[IM:23]] 
for $L^2$-subcritical problems
(e.g. $g(s) = \abs{s}^{p-1}s$, $1 < p < 1 + {4\over N}$) and
Bartsch-de Valeriola [\cite[BdV:3]], Bartsch-Soave [\cite[BS:5]], Bieganowski-Mederski [\cite[BM:8]],
Jeanjean [\cite[J:24]], Jeanjean-Lu [\cite[JL1:26]], Jeanjean-Le[\cite[JLe:25]],
Bartsch-Molle-Rizzi-Verzini [\cite[BMRV:4]] for $L^2$-supercritical problems
(e.g. $g(s) = \abs{s}^{p-1}s$, $1 + {4\over N} < p < 2^*-1$).
Here $2^* = {2N\over N-2}$ for $N \geq 3$ and $2^* = \infty$ for $N=2$.
See also [\cite[CGIT:15]] for a study of $L^2$-critical problems (c.f. Schino [\cite[Sc:36]],
Jeanjean-Zhang-Zhong [\cite[JZZ:31]]).
See also Jeanjean-Zhang-Zhong [\cite[JZZ:31]] for an approach using bifurcation.
We also refer to Bieganowski-d'Avenia-Schino [\cite[BDS:7]], 
Jeanjean-Lu [\cite[JL3:28]], Qi-Zou [\cite[QZ:34]], Soave [\cite[So:39]], Wei-Wu [\cite[WW:42]]
for studies of combined problems, that is, $g(s)$ is $L^2$-subcritical (resp. $L^2$-supercritical)
at $s=0$ and $L^2$-supercritical (resp. $L^2$-subcritical) at $s=\infty$.
We note that the orbital stability of solutions is also studied
in [\cite[CKS:10], \cite[CS:11], \cite[CL:12], \cite[Sh:38], \cite[TVZ:41]].


In this paper, we take another approach to $(*)_m$ and
we search for solutions of $(*)_m$ through a Lagrange
formulation of $(*)_m$:
    $$  \eqalignno{
        J_m(\mu,u) &= \half  \int_{\R^N} \abs{\nabla u}^2\, dx
             - \int_{\R^N} G(u)\, dx
             + \mu \left( \half  \int_{\R^N} u^2\, dx - m \right) \cr
        &:(0,\infty) \times H^1_r(\R^N) \to \R. &\label[1.2]\cr}
    $$
\noindent
Such approaches are successfully applied in 
Hirata-Tanaka [\cite[HT:22]] and Cingolani-Gallo-Ikoma-Tanaka [\cite[CGIT:15]] 
for $L^2$-subcritical and $L^2$-critical problems. 
We note that the least energy level, that is 
$\inf \{ \calI(u); u \in \calS_m \}$,
is characterized as a Mountain Pass (MP, in short) level for $J_m(\mu, u)$.
The Lagrangian approach is also successfully
applied to nonlinear Choquard problems in
Cingolani-Tanaka [\cite[CT1:16]], Cingolani-Gallo-Tanaka 
[\cite[CGT1:13], \cite[CGT3:14]].

Here we present new applications of our approach to nonlinear Schr\"odinger equations with 
cubic-quintic type nonlinearity and to $L^2$-supercritical problems.

\medskip

\BSS{\label[Subsection:1.1]. Applications}
We study the cubic-quintic problem and $L^2$-supercritical problems. 
In particular, when $N = 2,3$, we consider these problems under 
the following condition (p0) and we present the following existence results 
as special cases of our results in Section \ref[Section:5]

\smallskip

{\parindent=1.5\parindent
\item{(p0)} for some $p_0 \in (1 + {4\over N}, {N\over N-2}]$ when $N = 3$,  
$p_0 \in (1 + {4\over N},\, \infty)$ when $N = 2$,
    $$ \lim_{s \to 0^+} {g(s)\over s^{p_0}} = 1. 
    $$
}

\smallskip

\proclaim Theorem \label[Theorem:1.1]. (Cubic-quintic type problem)
Suppose $N = 2, 3$ and assume that $g(s)\in \m C([0,\infty),\R)$ satisfies (p0).
Moreover assume one of the following two conditions:
{
\item{(1)} $\displaystyle \limsup_{s\to+\infty} {g(s)\over s} \in [-\infty, \infty)$;
\item{or}
\item{(2)}  for some $p_\infty \in (1, 1 + {4\over N})$,  
    $$ \lim_{s\to+\infty} {g(s)\over s^{p_\infty}} = 1.
    $$
}
Then there exists $c_* \in (0,\infty)$ such that:
\item{(i)} if $m = c_*$, then $(*)_m$ has at least one positive solution  
$(\mu_1, u_1) \in (0,\infty) \times H^1_r(\R^N)$  
with  $J_m(\mu_1, u_1) > 0$.
\item{(ii)} if $m > c_*$, then $(*)_m$ has at least two positive solutions  
$(\mu_1,u_1), (\mu_2,u_2) \in (0,\infty) \times H^1_r(\R^N)$  
with  $ J_m(\mu_1,u_1) > 0$, $J_m(\mu_2,u_2) < J_m(\mu_1,u_1)$.

\medskip

We note that $1+{4\over N}<{N\over N-2}$ holds if and only if $N=2,3$.  We also note
that a characterization of $c_*$ in Theorem \ref[Theorem:1.1] and other existence
and non-existence results for $N\geq 4$ or $m<c_*$ will be given 
in Section \ref[Section:5].

\medskip

\proclaim Theorem \label[Theorem:1.2]. ($L^2$-supercritical problem)  
Suppose $N = 2,3$ and 
assume that $g(s)\in C([0,\infty),\R)$ satisfies (p0) and 
for $p_\infty \in (1 + {4\over N},\, 2^* - 1)$,
    $$ \lim_{s \to +\infty} {g(s)\over s^{p_\infty}} = 1.
    $$
Then for any $m > 0$, $(*)_m$ has at least one positive solution  
$(\mu_0, u_0) \in (0,\infty) \times H^1_r(\R^N)$  
with  $J_m(\mu_0, u_0) > 0$.

\medskip

In above two theorems, we assume that  
(a) $g(s)$ is $L^2$-supercritical at $s = 0$,  
(b) $g(s)$ is $L^2$-subcritical at $s = \infty$ (Theorem \ref[Theorem:1.1])  
or $L^2$-supercritical at $s = \infty$ (Theorem \ref[Theorem:1.2]).
We note that the assumption (1) of Theorem \ref[Theorem:1.1] holds for the cubic-quintic model $g(s)=s^3-s^5$, 
which appears in numerous problems in physics; field theory, nonlinear optics, vortex soliton theory,
Bose-Einstein condensates etc. (see [\cite[B:2], \cite[CKS:10], \cite[CS:11], \cite[KOPV:32]]).
In our theorems,  we assume only the conditions on the behaviors 
at $s\sim 0$, $s\sim\infty$ and the continuity of $g(s)$ in $(0,\infty)$
and we do not assume any other condition on $g(s)$.
In particular our Theorem \ref[Theorem:1.2] accepts sign-changing $g(s)$.
See Section \ref[Section:5] for other existence results for $L^2$-supercritical problems.

In the preceding results [\cite[JL3:28]] for the cubic-quintic problem and 
[\cite[BS:5], \cite[BM:8], \cite[J:24], \cite[JL1:26]]
they require some global conditions,
typically the global Ambrosetti-Rabinowitz type conditions:
    $$  \theta G(s) \ge g(s) s \quad \hbox{for all}\  s \in \R, 
    $$
for some $\theta \in (2, 2^*]$. 
We do not require such global conditions and thus
our Theorems \ref[Theorem:1.1] and \ref[Theorem:1.2] improve the preceding results when $N=2,3$.
See Section \ref[Section:5] for the existence results for $N\geq 4$.

Our Theorems \ref[Theorem:1.1] and \ref[Theorem:1.2] will be obtained from the general theory 
for $J_m(\mu, u)$.  

\medskip


\BSS{\label[Subsection:1.2]. General theory}
We develop an abstract theory for $J_m(\mu,u)$ in a general situation and we assume
the following conditions

\smallskip

{\parindent=1.5\parindent

\item{(g0)} $g(s) \in C([0,\infty), \R)$;

\item{(g1)} $\lim_{s \to 0^+} {g(s)\over s} = 0$;

\item{(g2)} when $N \geq 3$,
    $\displaystyle  \lim_{s \to \infty} {g(s)\over s^{2^*-1}} = 0
    $, \m
when $N = 2$, for any $\alpha > 0$,
    $\displaystyle  \lim_{s \to \infty} {g(s)\over e^{\alpha s^2}} = 0
    $;

\item{(g3)} there exists $s_0 > 0$ such that $G(s_0) > 0$.

} 
\smallskip

\noindent
We assume conditions (g0)--(g3) throughout this paper.
For technical reasons we set
    $$  g(s) = -s \quad \hbox{for } s \in (-\infty,0)   \eqno\label[1.3]
    $$
and we look for positive solutions.  

For $\mu > 0$, we consider the nonlinear scalar field
equation (with a fixed frequency)
    $$  (**) \quad -\Delta u + \mu u = g(u) \quad \hbox{in } \R^N.
    $$
Nonlinear scalar fields are well-studied since the
pioneering work Berestycki-Lions [\cite[BL:6]] and under conditions closely
related to (g0)--(g3) the existence of a least energy
solution is shown. 
Moreover, in Jeanjean-Tanaka [\cite[JT1:29], \cite[JT2:30]] it is
shown that the least energy solutions enjoy
a MP minimax characterization using
the free functional. In addition, using a new
version of the Palais-Smale condition $(PSP)$
(see Sections \ref[Section:2], \ref[Section:3]), the attainability of the MP 
level is shown in [\cite[HT:22]], Cingolani-Tanaka [\cite[CT2:17]] 
(cf. Hirata-Ikoma-Tanaka [\cite[HIT:21]]).

In our setting we can show that $(**)$ satisfies
the conditions in [\cite[BL:6]] for $\mu \in (0,\mu_*)$, where
$\mu_*$ is defined by
    $$  \mu_* = \sup_{s>0} {G(s)\over {1\over 2}s^2} \in (0,\infty] ,
        \eqno \label[1.4]
    $$
and thus $(**)$ has a least energy solution for $\mu \in (0,\mu_*)$.
We note that $\mu_*<\infty$ if $g(s)$ only has a linear growth rate from above,
that is,  $\limsup_{s\to\infty}{g(s)\over s}<\infty$.

Solutions of $(**)$ are characterized as critical points
of the following free functional:
    $$  \eqalignno{
        I(\mu,u) &= \half \int_{\R^N} \abs{\nabla u}^2 \, dx
        + {\mu\over 2} \int_{\R^N} u^2 \, dx
        - \int_{\R^N} G(u)\, dx \cr
        &: (0,\infty) \times H^1_r(\R^N) \to \R. &\label[1.5]\cr}
    $$
Thus for $\mu \in (0,\mu_*)$ the least energy level is
characterized as a MP level
    $$  a(\mu) = \inf_{\gamma \in \Gamma_\mu} \max_{\tau \in [0,1]} I(\mu, \gamma(\tau)),
        \eqno\label[1.6]
    $$
where
    $$  \Gamma_\mu = \{ \gamma \in C([0,1], H^1_r(\R^N)) ;
        \ \gamma(0) = 0,\ I(\mu,\gamma(1)) < 0 \}.
    $$
We will see in Section \ref[Subsection:2.2] that 
$a(\mu)$ is attained by a critical point of 
$u\mapsto I(\mu,u)$ and the set of critical points
    $$  \calC(\mu) = \{ u \in H^1_r(\R^N) ; I(\mu,u) = a(\mu), \, \partial_u I(\mu,u)=0 \} 
        \eqno\label[1.7]
    $$
is compact and non-empty.
See Section \ref[Subsection:2.2] for more properties of $a(\mu)$.  

We note $J_m(\mu,u) = I(\mu,u) - m\mu$. Thus
    $$  b_m(\mu) = a(\mu) - m\mu    \eqno\label[1.8]
    $$
is the least energy level for $u \mapsto J_m(\mu,u)$. 
In particular, 
    $$  J_m(\mu,u)\geq b_m(\mu), \quad \hbox{i.e.,}\quad \calI(u)\geq b_m(\mu)
        \quad \hbox{for all solutions $(\mu,u)$ of $(*)_m$.}
                \eqno\label[1.9]
    $$
Through $b_m(\mu)$, we will show the existence of critical points of $J_m(\mu,u)$, 
i.e., solutions of $(*)_m$.

\medskip

First we study the differentiability of $a(\mu)$. We need the following notation 
for $\mu \in (0,\mu_*)$:
    $$  c_+(\mu) = \max_{u \in \calC(\mu)} \half  \int_{\R^N} u^2 dx, \quad
        c_-(\mu) = \min_{u \in \calC(\mu)} \half  \int_{\R^N} u^2 dx, 
        \eqno\label[1.10]
    $$
where $\calC(\mu)$ is the compact set defined in \ref[1.7].  
We have the following result.

\medskip

\proclaim Theorem \label[Theorem:1.3].  
Assume (g0)--(g3) and \ref[1.3]. Let $a(\mu)$, $c_+(\mu)$, $c_-(\mu)$ be functions
defined in \ref[1.6], \ref[1.10].  
Then $a(\mu)$ is right and left differentiable and
    $$  a'_+(\mu) = c_-(\mu), \quad a'_-(\mu) = c_+(\mu) 
        \qquad \hbox{for } \mu \in (0,\mu_*).
    $$

\medskip

As a corollary to Theorem \ref[Theorem:1.3], we have

\medskip

\proclaim Theorem \label[Theorem:1.4].
Assume (g0)--(g3) and \ref[1.3].  
For $m>0$, 
let $b_m(\mu)$ be given in \ref[1.8].  
If $\mu_0 \in (0,\mu_*)$ is a local minimizer of $b_m(\mu)$, then
\item{(i)} $b_m(\mu)$ is differentiable at $\mu=\mu_0$ and 
$(b_m)_+'(\mu_0)=(b_m)_-'(\mu_0)=0$;
\item{(ii)} for any $u_0\in\calC(\mu_0)$, $(\mu_0,u_0)$ is 
a positive solution of $(*)_m$ and
    $$  I(\mu_0,u_0) = b_m(\mu_0).
    $$

\medskip


\claim Remark \label[Remark:1.5]. {\sl
Let $\mu \in (0,\mu_*)$ and $u \in \calC(\mu)$.
If $(\mu,u)$ solves $(*)_m$, we have
    $$  c_{-}(\mu) \le \half  \norm{u}_2^2 = m \le c_{+}(\mu).
    $$
Noting
    $   (b_m)'_{\pm}(\mu) = c_{\mp}(\mu) - m 
    $,
we have $(b_m)'_{+}(\mu) \le (b_m)'_{-}(\mu)$.
Thus

\item{(i)} if $(b_m)'_{+}(\mu) > 0$
(resp. $(b_m)'_{-}(\mu) < 0$),
then $c_{-}(\mu) > m$ (resp. $c_{+}(\mu) < m$) and
    $$  \half  \norm{u}_2^2 > m 
        \quad (\hbox{resp.} <m) \quad \hbox{for all}\  u \in \calC(\mu).
    $$
Thus $(\mu,u)$ is not a solution of $(*)_m$
for all $u \in \calC(\mu)$;

\item{(ii)} if $(b_m)'_{+}(\mu)=0$ (resp. $(b_m)'_{-}(\mu)=0$),
clearly for $u \in \calC(\mu)$ with $\half\norm u_2^2 = c_-(\mu)$ 
(resp. $\half\norm u_2^2 = c_+(\mu)$),
$(\mu,u)$ solves $(*)_m$. 
Moreover, if $(b_m)'_{+}(\mu) = (b_m)'_{-}(\mu)=0$,
then for any $u \in \calC(\mu)$,
$(\mu,u)$ solves $(*)_m$.  We note that 
$(b_m)'_{+}(\mu) = (b_m)'_{-}(\mu)=0$ holds at a local minimum of $b_m$;

\item{(iii)} if $(b_m)'_{+}(\mu) < 0 < (b_m)'_{-}(\mu)$,
we observe that $\mu$ is a local maximum of $b_m$.
However, the existence of $u \in \calC(\mu)$
with $\half  \norm{u}_2^2 = m$
is not clear.
Our Theorem \ref[Theorem:1.6] below deals with this case
under additional assumptions.

} 


\medskip

It is natural to expect the existence of a critical point of $J_m(\mu,u)$
which is related to a local maximum of $b_m(\mu)$.  
The following theorem gives a partial affirmative answer.

\medskip

\proclaim Theorem \label[Theorem:1.6].  
Assume (g0)--(g3), \ref[1.3] and let $b_m(\mu)$ be given in \ref[1.8].  
Assume that $b_m(\mu)$ takes a topologically non-trivial local maximum 
at $\mu_0 \in (0,\mu_*)$, that is, for some interval 
$[\mu_1,\mu_2] \subset (0,\mu_*)$ we assume $\mu_0 \in (\mu_1,\mu_2)$ and
    $$  b_m(\mu_0) = \max_{\mu \in [\mu_1,\mu_2]} b_m(\mu)
        > \max\{ b_m(\mu_1),\ b_m(\mu_2)\}.
    $$
Moreover, assume also $\widetilde J_m(\lambda,u) = J_m(e^\lambda,u)$ satisfies 
the $(PSPC)_b$ condition given in Definition \ref[Definition:2.8] 
in Section \ref[Subsection:2.3] for all $b \geq b_m(\mu_0)$.  
Then $J_m(\mu,u)$ has a critical point $(\mu_\#, u_\#)$, i.e., 
a positive solution of $(*)_m$, such that
    $$  J_m(\mu_\#,u_\#) \geq b_m(\mu_0).
    $$

\medskip

We use ideas in Cingolani-Gallo-Ikoma-Tanaka [\cite[CGIT:15]] to show Theorem \ref[Theorem:1.6].  
We will see that a topologically non-trivial local maximum of $b_m(\mu)$ gives 
a local linking structure of $\widetilde J_m(\lambda,u)$ in $\R\times H^1_r(\R^N)$
and $b_m(\mu)$ gives a new way to find a local linking 
for $J_m(\mu, u)$. 
The $(PSPC)_b$ condition enables us to use a deformation argument for 
$\widetilde J_m(\lambda,u)$ in $\R\times H^1_r(\R^N)$.

Our Theorems \ref[Theorem:1.1] and \ref[Theorem:1.2] (more generally Theorems \ref[Theorem:5.1] and \ref[Theorem:5.4]) 
are consequences of our abstract results (Theorems \ref[Theorem:1.4] and \ref[Theorem:1.6]) and studies of the behaviors 
of $b_m(\mu)$, which will be given in Sections \ref[Section:2]--\ref[Section:4].  
The details of the proofs will be given in Section \ref[Section:5].
Here we give typical figures of $b_m(\mu)$ in the setting of Theorems \ref[Theorem:1.1] and \ref[Theorem:1.2].

\bigskip
\medskip

\line{
\hfil\pdffig{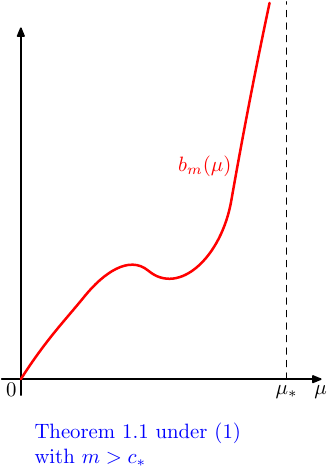}\hfil
\hfil\pdffig{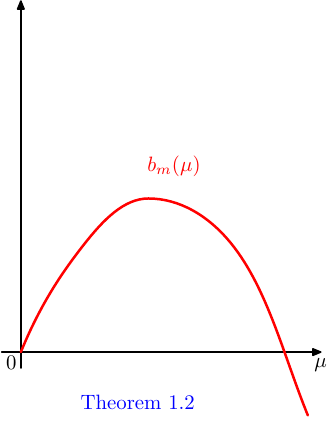}
\hfil}

\medskip


\medskip

\claim Remark \label[Remark:1.7].
{\sl 
The present paper is largely motivated by 
Dovetta-Serra-Tilli [\cite[DST2:20]], in which they 
study the following prescribed mass problem for
    $$  \left\{ \eqalign{
        -&\Delta u + \mu u = \abs{u}^{p-1}u \quad \hbox{in } \Omega, \cr
        &u=0 \qquad  \hbox{on } \partial \Omega, \cr}
        \right.
        \eqno\label[1.11]
    $$
where $\Omega \subset \R^N$ is a (possibly unbounded) open set and $p \in (1,2^*)$. 
They introduce two important energy levels, action and energy ground state levels, 
and they study properties and relations of these two energy levels. 

Their action ground state level is introduced as
    $$  \widetilde a(\mu) = \inf \{ \widetilde I(\mu,u) \,|\, u \in \widetilde\calN_\mu \}, 
    $$
where the functional $\widetilde I$ and the Nehari manifold $\widetilde \calN_\mu$
are associated to \ref[1.11].  Namely,
    $$  \eqalign{
        &\widetilde I(\mu,u) = {1\over 2}\int_\Omega \abs{\nabla u}^2 dx 
            + {\mu\over 2}\int_\Omega \abs{u}^2 dx 
            - {1\over p+1}\int_\Omega \abs{u}^{p+1} dx, \cr
        &\widetilde\calN_\mu = \{ u \in H^1_0(\Omega)\setminus\{0\} ;
                    \, \partial_u \widetilde I(\mu,u)u = 0 \}. \cr}
    $$
We note that $\widetilde a(\mu)$ is nothing but the MP level for $u \mapsto I(\mu,u)$.

Dovetta-Serra-Tilli [\cite[DST2:20]] (c.f. [\cite[DST1:19]]) 
investigate the differentiability of $\widetilde a(\mu)$, which is closely related 
to our Theorem \ref[Theorem:1.3]. 
They also obtain the existence of a prescribed mass solution corresponding to 
a minimum and a local minimum of 
    $   \widetilde b_m(\mu) = \widetilde a(\mu) - m\mu
    $.

We have shown in Theorem \ref[Theorem:1.3]
the existence and representation 
of the right and left derivatives of $a(\mu)$ for general nonlinearity $g(s)$, 
which enables us 
to study the existence of critical points related to local minima of 
$b_m(\mu)$ (Theorem \ref[Theorem:1.4]). 
Moreover, we study the existence of critical points related to local maxima 
of $b_m(\mu)$ (Theorem \ref[Theorem:1.6]).

} 

\medskip

This paper is organized as follows.  In Section \ref[Section:2], we give some preliminary 
results on prescribed frequency problems and Palais-Smale type conditions.
To apply our abstract results (Theorems \ref[Theorem:1.4] and \ref[Theorem:1.6]) to $(*)_m$,
the properties of $b_m(\mu)$ and the $(PSPC)$ condition are important.
In Section \ref[Section:3], we study the behaviors of the MP value $a(\mu)$
and $c_-(\mu)$ as $\mu\sim 0$ and $\mu\sim \infty$.
In Section \ref[Section:4], we study $(PSPC)$ condition.  In Sections \ref[Section:3], \ref[Section:4],
the behavior of $g(s)$ at $s\sim 0$ and $s\sim \infty$ are essential.
In Section \ref[Section:5], we give our existence results for $(*)_m$ --- special cases of
our existence results are given in the Introduction as Theorems \ref[Theorem:1.1] and \ref[Theorem:1.2] ---
using the results in Sections \ref[Section:3], \ref[Section:4] and our abstract results.
Finally in Section \ref[Section:6], we prove our abstract results 
(Theorems \ref[Theorem:1.3], \ref[Theorem:1.4] and \ref[Theorem:1.6]).

\medskip

\BS{\label[Section:2]. Preliminaries}
Throughout this paper, we assume (g0)--(g3) and \ref[1.3].
We start with the variational formulation of our problem.

\BSS{\label[Subsection:2.1] Variational formulation and function spaces}
We look for radially symmetric solutions of $(*)_m$ and
work in the space of radially symmetric functions:
    $$  E = H^1_r(\R^N) 
        = \{ u \in H^1(\R^N);\  u(x) = u(\abs{x}) \}.
    $$
We use the following notation for $u$, $v \in E$:
    $$  \eqalign{
        &\norm{u}_p = \Bigl( \int_{\R^N} \abs{u(x)}^p dx \Bigr)^{1/p}
            \quad \hbox{for}\ p \in [1,\infty), \qquad
        \norm{u}_\infty = \esssup_{x \in \R^N} \abs{u(x)}, \cr
        &\norm{u}_E = \bigl( \norm{\nabla u}_2^2 + \norm{u}_2^2 \bigr)^{1/2}, \quad
        (u,v)_2 = \int_{\R^N} u v \, dx. \cr}
    $$
We also use the notation
    $$  u_+(x)=\max\{ u(x), 0\}, \quad u_-(x)=\max\{ -u(x), 0\}
    $$
for positive and negative parts of $u$.

For a given $m>0$, we take a variational approach
to $(*)_m$ and we search for a critical point of the following functional:
    $$  J_m(\mu,u) = \half \norm{\nabla u}_2^2 + {\mu\over 2}\norm{u}_2^2 
        + \mu \Bigl(\half \norm{u}_2^2 - m\Bigr)
        \in C^1((0,\infty)\times E,\R).
    $$
We note that
    $$  \eqalignno{
        &\partial_\mu J_m(\mu,u) = \half \norm{u}_2^2 - m, \cr
        &\partial_u J_m(\mu,u)\varphi 
            = (\nabla u, \nabla \varphi)_2 + \mu (u,\varphi)_2 
            - \int_{\R^N} g(u)\varphi \, dx
            \quad \hbox{for all}\  \varphi \in E.       &\label[2.1] \cr}
    $$
Thus $(\mu,u)\in(0,\infty)\times E$ is a critical point of $J_m(\mu,u)$
if and only if $(\mu,u)$ is a solution of $(*)_m$.

Setting $\varphi(x) = u_-(x)$ in \ref[2.1], we have 
    $   \norm{\nabla u_-}_2^2 + (\mu+1)\norm{u_-}_2^2 = 0
    $
from the assumption \ref[1.3].
Thus $u_-\equiv 0$ and $u(x)$ is a non-negative solution.
Since $u(x)$ solves $-\Delta u + \mu u = g(u)$ in $\R^N$, 
we can also observe that $u(x)>0$ in $\R^N$ if $u(x)$ is not identically zero.

\medskip

It is also well-known that solutions $(\mu,u)\in(0,\infty)\times E$
satisfy the Pohozaev identity:
    $$  P(\mu,u)=0,
    $$
where $P(\mu,u)\in C^1((0,\infty)\times E,\R)$ is defined by
    $$  P(\mu,u) = {N-2\over 2}\norm{\nabla u}_2^2 
        + N\Bigl({\mu\over 2}\norm{u}_2^2 - \int_{\R^N} G(u)dx\Bigr).
    $$
Formally we have
    $$  P(\mu,u) = \left.{d\over d\tau}\right|_{\tau=1} J_m(\mu, u(x/\tau))
    $$
and functional $P(\mu,u)$ plays important roles in our deformation argument.


\BSScap{\label[Subsection:2.2]. Problems with prescribed frequency and MP value $a(\mu)$}{\ref[Subsection:2.2]. Problems with prescribed frequency and MP value a}
To study $(*)_m$, first we consider the prescribed frequency problem,  
that is, for fixed $\mu > 0$ we consider the problem $(**)$.  
The functional corresponding to $(**)$ is
    $$  u \mapsto I(\mu,u) = \half\norm{ \nabla u }_2^2 
            - \int_{\R^N} G_\mu(u)\, dx ,
    $$
where $G_\mu(s) = G(s) - {\mu\over 2} s^2$.  

Defining $\mu_* \in (0,\infty]$ by \ref[1.4], i.e.,
$\mu_* = \sup_{s>0} {G(s)\over \half s^2} \in (0,\infty]$, we have
the following facts.

\smallskip

\item{(1)} If $\mu_* < \infty$,  
    $$  G(s) \leq {\mu_*\over 2} s^2 \quad \hbox{for all}\  s \geq 0.
        \eqno\label[2.2]
    $$
\item{(2)} $G_\mu(s)$ satisfies the conditions in [\cite[BL:6]],  
i.e., $\sup_{s>0} G_\mu(s) > 0$, if and only if $\mu \in (0,\mu_*)$.  
\item{(3)} For $\mu \in (0,\mu_*)$, $u \mapsto I(\mu,u)$ has a MP geometry, that is,  
\itemitem{(i)} $I(\mu,0) = 0$;  
\itemitem{(ii)} there exist $r_\mu>0$ and $c_\mu > 0$ such that 
$I(\mu,u) \geq c_\mu$ for $u \in E$ with $\norm{u}_E = r_\mu$;
\itemitem{(iii)} there exists $e_\mu \in E$ such that $\norm{e_\mu}_E > r_\mu$ 
and $I(\mu,e_\mu) < 0$.

\smallskip

\noindent
By (3), we can define a MP value for
$u \mapsto I(\mu,u)$ for each $\mu \in (0,\mu_*)$ by
    $$  \eqalignno{
        &a(\mu) = \inf_{\gamma \in \Gamma_\mu} \max_{\tau \in [0,1]} 
            I(\mu, \gamma(\tau)) (\geq c_\mu>0),
            &\label[2.3]\cr
        &\Gamma_\mu = \{ \gamma(\tau) \in C([0,1],E) ; 
                \ \gamma(0)=0,\ \gamma(1)=e_\mu \}. &\label[2.4]\cr}
    $$
We note that the set $\{ u \in E ; \ I(\mu,u) < 0 \}$ is path-connected 
(c.f. [\cite[By:9], \cite[HIT:21]]).  
The MP value $a(\mu)$ does not depend on the choice of the end point $e_\mu$.  

Recalling the results in [\cite[CT2:17], \cite[HIT:21], \cite[JT1:29]], we have

\medskip

\proclaim Lemma \label[Lemma:2.1].
For each $\mu \in (0,\mu_*)$, we have  
\item{(i)} $a(\mu)$ is attained by a critical point;
\item{(ii)} $a(\mu)$ is the least energy level for $u \mapsto I(\mu,u)$, that is,  
    $$  a(\mu) = \inf \{ I(\mu,u); \, u \neq 0, \ \partial_u I(\mu,u) = 0 \};
    $$  
\item{(iii)} the critical set $\calC(\mu) 
= \{ u \in E; \, \partial_u I(\mu,u)=0,\ I(\mu,u)=a(\mu) \}$ is compact and non-empty.
Moreover 
    $$  \bigcup_{\mu \in [\mu_1,\mu_2]} \left( \{ \mu \} \times \calC(\mu) \right)
    $$  
is compact in $(0,\infty) \times E$ for any compact interval 
$[\mu_1,\mu_2] \subset (0,\mu_*)$;
\item{(iv)} $(0,\mu_*) \to (0,\infty);\ \mu \mapsto a(\mu)$ is continuous and 
strictly increasing;
\item{(v)} $c_-(\mu) = \min_{u \in C(\mu)} \half  \norm{u}_2^2$,  
    $c_+(\mu) = \max_{u \in C(\mu)} \half  \norm{u}_2^2$  
are well-defined and positive for each $\mu \in (0,\mu_*)$.

\medskip

\claim Proof.
(i)--(iii) are essentially given in [\cite[JT1:29]].  
The attainability of $a(\mu)$ is also proved in [\cite[HT:22], Section \raw[6]], [\cite[CT2:17]].
(iv)--(v) also follow easily.
\QED

\medskip

We will study the behavior of $a(\mu)$ as $\mu\sim 0$ and $\mu\sim\mu_*$ precisely 
in Section \ref[Section:3].

Finally in this section, we recall another characterization of 
the least energy level $a(\mu)$.  

\medskip

\proclaim Proposition \label[Proposition:2.2]. 
For $\mu\in (0,\mu_*)$,
    $$  a(\mu) = \inf_{u\in \calP(\mu)} I(\mu,u) 
        = \inf_{u\in \calP(\mu)} {1\over N}\norm{\nabla u}_2^2 >0,
    $$
where $\calP(\mu)$ is the Pohozaev manifold for $-\Delta u+\mu u=g(u)$ defined
by
    $$  \calP(\mu) = \{u\in E\setminus \{0\};\ P(\mu,u)=0\}. 
        \eqno \label[2.5]
    $$

\medskip

\claim Proof.
We note that
    $$  I(\mu,u) = {1\over N}\norm{\nabla u}_2^2 \quad \hbox{on}\ \calP_\mu.
    $$
Proposition \ref[Proposition:2.2] is a special case of the results in [\cite[BL:6], \cite[CGM:18], \cite[JT1:29]].  
\QED

\medskip


\BSS{\label[Subsection:2.3]. Palais-Smale type conditions and deformation results}
We use deformation arguments to study $(**)$ and $(*)_m$ as 
in [\cite[HT:22], Section \raw[6]], [\cite[CT2:17]] and [\cite[CGIT:15]].  
In these works the scaling properties of the functionals 
are very important and new versions of Palais-Smale condition play essential roles.

\medskip

\noindent
{\bf a) Palais-Smale type condition and deformation argument for $(**)$}

\medskip

\noindent
We give a review of results in [\cite[CT2:17], \cite[HT:22]].
We consider
    $$  - \Delta u = \widehat{g}(u) \quad \hbox{in}\  \R^N \eqno \label[2.6]
    $$
For $\widehat{g}(s)$, we assume (g0), (g2), (g3) and  
    $$  -\infty < \lim_{s \to 0^+} {\widehat{g}(s)\over s} 
        < 0.
        \eqno\label[2.7]
    $$
We also assume that $\widehat{g}(s)$ satisfies \ref[1.3].

The functional corresponding to \ref[2.6] is
    $$  F(u) = \half  \norm{ \nabla u }_2^2 - \int_{\R^N} \widehat{G}(u)\, dx \in C^1(E,\R),
    $$
where $\widehat{G}(s) = \int_0^s \widehat{g}(\tau)\, d\tau$.

\medskip

\claim Definition \label[Definition:2.3] {\rm ($(PSP)$ condition in $E$)}.  
{\sl
For $b \in \R$, we say that $F(u)$ satisfies the {\it Palais-Smale-Pohozaev condition}
at level $b$ in $E$ ($(PSP)_b$ condition, in short), if any sequence 
$(u_j)_{j=1}^\infty \subset E$ satisfying as $j \to \infty$,
    $$  \eqalignno{
        &F(u_j) \to b, &\label[2.8]\cr
        &\norm{ \partial_u F(u_j) }_{E^*} \to 0, &\label[2.9]\cr
        &Q(u_j) \to 0 &\label[2.10]\cr}
    $$
has a strongly convergent subsequence in $E$.  
Here $Q(u) \in C^1(E,\R)$ is the Pohozaev functional associated to \ref[2.6], i.e.,
    $$  Q(u) = {N-2\over 2} \norm{ \nabla u }_2^2 - N \int_{\R^N} \widehat{G}(u)\, dx .
    $$
We also say that $(u_j)_{j=1}^\infty \subset E$ is a {\it $(PSP)_b$ sequence} for $F(u)$, 
if $(u_j)_{j=1}^\infty$ satisfies \ref[2.8]--\ref[2.10].
} 

We developed the following deformation result in [\cite[HT:22], Section \raw[6]], [\cite[CT2:17]].

\medskip

\proclaim Proposition \label[Proposition:2.4].
{\rm ([\cite[CT2:17], Theorem \raw[2.1]])}
For $b \in \R$, assume that $F(u)$ satisfies the $(PSP)_b$ condition. 
Then for any neighborhood $\calO$ of the critical set  
$K_b = \{ u \in E; F(u)=b, \ \partial_u F(u)=0,\ Q(u)=0 \}$ 
and for any $\bar{\epsilon} > 0$ there exist $\epsilon \in (0,\bar{\epsilon})$ 
and a continuous map $\eta(t,u) : [0,1] \times E \to E$ such that  
\item{(i)} $\eta(0,u) = u$ for all $u \in E$;  
\item{(ii)} $\eta(t,u) = u$ for all $t \in [0,1]$ if $F(u) \leq b-\bar\epsilon$;  
\item{(iii)} $F(\eta(t,u)) \leq F(u)$ for all $(t,u) \in [0,1] \times E$;  
\item{(iv)} $F(\eta(1,u)) \leq b-\epsilon$, if $u \in E \setminus \calO$ and 
$F(u) \leq b+\epsilon$;  
\item{(v)} if $K_b = \emptyset$, then $F(\eta(1,u)) \leq b-\epsilon$ 
for $u \in E$ with $F(u) \leq b+\epsilon$.

\medskip

We also have

\medskip

\proclaim Proposition \label[Proposition:2.5]. {\rm ([\cite[HT:22], Proposition \raw[6.1]],
c.f. [\cite[HIT:21]])}
Assume that $\widehat{g}(s)$ satisfies (g0), (g2), (g3), \ref[2.7] and \ref[1.3].  
Then for any $b \in \R$, $F(u)$ satisfies the $(PSP)_b$ condition.

\medskip

\claim Remark \label[Remark:2.6].
{\sl 
In [\cite[CT2:17], \cite[HT:22]], it is assumed that $\widehat{g}(s)$ is odd and 
results corresponding to Propositions \ref[Proposition:2.4] and 
\ref[Proposition:2.5] are shown. 
After slight modifications, the arguments work under the assumption \ref[1.3].

\medskip

For any $\mu \in (0,\mu_*)$, Propositions \ref[Proposition:2.4] and \ref[Proposition:2.5] 
are applicable to 
$F(u) = I(\mu,u)$. Thus we can show the MP minimax value $a(\mu)$ defined 
in \ref[2.3]--\ref[2.4] is a critical value of $u \mapsto I(\mu,u)$.
} 

\medskip

Modifying the proof of Proposition \ref[Proposition:2.5] slightly, we can show the following.

\medskip

\proclaim Proposition \label[Proposition:2.7].
Let $[\mu_1,\mu_2] \subset (0,\infty)$ be a compact interval.  
Assume $(\mu_j, u_j)_{j=1}^\infty \subset (0,\mu_*) \times E$ satisfy  
    $$  \eqalign{
        &\mu_j \in [\mu_1,\mu_2] \quad \hbox{for all}\  j \in \N,  \cr
        &I(\mu_j,u_j) \to b, \cr
        &\norm{\partial_u I(\mu_j,u_j)}_{E^*} \to 0,  \cr
        &P(\mu_j,u_j) \to 0 \cr}
    $$
as $j\to\infty$.
Then $(\mu_j,u_j)_{j=1}^\infty$ has a strongly convergent subsequence 
in $(0,\infty) \times E$.


This proposition is useful to verify $(PSPC)$ condition
for $J_m(\mu, u)$ in Section \ref[Section:4].

\medskip

\noindent
{\bf b) Palais-Smale type condition and deformation
argument for $(*)_m$}

\medskip

\noindent
In the deformation argument, completeness of the base space is important. 
Introducing a change of variable $\mu = e^\lambda$,
we introduce the following functional defined
on a complete space $\R\times E$,
    $$  \widetilde{J}_m(\lambda, u) = J_m(e^\lambda, u) 
        = \half \norm{\nabla u}_2^2 - \int_{\R^N} G(u)\,dx 
            + e^\lambda \left(\half \norm{u}_2^2 - m\right)
        \in C^1(\R\times E, \R). 
    $$
Clearly, there is one-to-one correspondence between critical points of $J_m$
and $\widetilde J_m$.
We also introduce
    $$  \widetilde{P}(\lambda, u) = P(e^\lambda, u) 
        = {N-2\over 2}\norm{\nabla u}_2^2 + N\left({e^\lambda\over 2}\norm{u}_2^2 
            - \int_{\R^N} G(u)\,dx\right) \in C^1(\R\times E, \R). 
    $$
The following Palais-Smale-Pohozaev-Cerami condition
at level $b$ in $\R\times E$ ($(PSPC)_b$ condition, in short)
is introduced in [\cite[CGIT:15]]. See also [\cite[HT:22]] for $(PSP)$
condition in $\R\times E$.

\medskip

\claim Definition \label[Definition:2.8] {\rm ($(PSPC)$ condition in $\R\times E$)}.
{\sl 
For $b\in\R$, we say that $\widetilde{J}_m(\lambda,u)$ satisfies the
{\it Palais-Smale-Pohozaev-Cerami condition} at level $b$ in $\RE$
({\it $(PSPC)_b$ condition} in $\R\times E$), if any sequence
$(\lambda_j, u_j)_{j=1}^\infty \subset \R\times E$ satisfying, as $j\to\infty$
    $$  \eqalignno{
        &\widetilde J_m(\lambda_j, u_j) \to b, &\label[2.11]\cr
        &\partial_\lambda \widetilde J_m(\lambda_j,u_j) \to 0, &\label[2.12]\cr
        &(1+\norm{u_j}_E)\norm{\partial_u \widetilde J_m(\lambda_j,u_j)}_{E^*} \to 0, 
                        &\label[2.13]\cr
        &\widetilde{P}(\lambda_j,u_j) \to 0 &\label[2.14]\cr}
    $$
has a strongly convergent subsequence in $\R\times E$.

We also say that $(\lambda_j,u_j)_{j=1}^\infty \subset \R\times E$ is a {\it $(PSPC)_b$
sequence} for $\widetilde J_m(\lambda,u)$ if $(\lambda_j,u_j)_{j=1}^\infty$ satisfies
\ref[2.11]--\ref[2.14].
} 

The following deformation result is essentially given in [\cite[CGIT:15], 
Proposition \raw[7.1]].  The situation in [\cite[CGIT:15]] is slightly different 
from this paper, but it is easy to modify.

\medskip

\proclaim Proposition \label[Proposition:2.9].  
{\rm ([\cite[CGIT:15], Proposition \raw[7.1]])}
Let $m>0$.  For $b\in\R$ assume that $\widetilde J_m(\lambda,u)$ satisfies
the $(PSPC)_b$ condition.  Then for any $\bar\epsilon > 0$ and
for any neighborhood $\calO$ of 
$\widetilde{K}_b = \{(\lambda,u)\in\R\times E;\, 
\widetilde J_m(\lambda,u)=b,\,\partial_\lambda \widetilde J_m(\lambda,u)=0,\,
\partial_u \widetilde J_m(\lambda,u)=0,\ \widetilde P(\lambda,u)=0\}$
in $\RE$, 
there exist $\epsilon\in(0,\bar\epsilon)$ and a continuous map
$\eta(t,\lambda,u): [0,1]\times\R\times E \to \R\times E$ such that
\item{(i)} $\eta(0,\lambda,u)=(\lambda,u)$ for all $(\lambda,u)\in\R\times E$;
\item{(ii)} 
$\eta(t,\lambda,u)=(\lambda,u)$ for all $t\in[0,1]$ if $\widetilde J_m(\lambda,u)\leq b-\bar\epsilon$;
\item{(iii)} $\widetilde J_m(\eta(t,\lambda,u)) \leq \widetilde J_m(\lambda,u)$ for all $t\in[0,1]$
and $(\lambda,u)\in\R\times E$;
\item{(iv)} $\widetilde J_m(\eta(1,\lambda,u)) \leq b-\epsilon$ if $(\lambda,u)\in \RE\setminus \calO$
and $\widetilde J_m(\lambda,u)\leq b+\epsilon$;
\item{(v)} if $\widetilde{K}_b=\emptyset$, then 
$\widetilde J_m(\eta(1,\lambda,u)) \leq b-\epsilon$
for $(\lambda,u)\in\R\times E$ with $\widetilde J_m(\lambda,u)\leq b+\epsilon$.

\medskip

We will use Proposition \ref[Proposition:2.9] to show Theorem \ref[Theorem:1.6].
In contrast to the prescribed frequency problem \ref[2.6],
verification of the $(PSPC)_b$ condition is a delicate problem and the behavior of
$g(s)$ at $0$ and $\infty$ is important.
We will study the $(PSPC)_b$ condition in Section \ref[Section:4].

\medskip

\claim Remark \label[Remark:2.10].  
{\sl When $g(s)=s^{1+{4\over N}}$, let
$\omega_1(x)\in E$ be the least energy solution of
$-\Delta u+u=u^{1+{4\over N}}$ and 
$m_1=\half \norm{\omega_1}_2^2$.
Then the corresponding functional
does not satisfy the $(PSPC)_b$ condition for $b=0$.
In fact, setting $u_\lambda(x)=\mu^{N/4}\omega_1(\mu^{1/2}x),\;\mu=e^\lambda$,
we observe that $(\lambda,u_\lambda)$ is a critical point of $\widetilde J_{m_1}$
at level $0$ for all $\lambda\in\R$. Thus $\widetilde J_{m_1}$ does not
satisfy the $(PSPC)_0$ condition.  See [\cite[CGIT:15]].

} 

\medskip

\BScap{\label[Section:3]. Behaviors of $a(\mu)$ and $c_{-}(\mu)$}{\ref[Section:3]. Behaviors of a and c}
To find solutions of $(*)_m$ through $b_m(\mu)$ using Theorems \ref[Theorem:1.4] 
and \ref[Theorem:1.6], we need to study the behavior of $b_{m}(\mu)$.  
The following propositions provide fundamental information on the behaviors of 
$a(\mu)$ and $c_{-}(\mu)$.

\medskip

\BSScap{\label[Subsection:3.1].  Behavior as $\mu \to \infty$, when $\mu_* = \infty$}{\ref[Subsection:3.1]. Behavior as mu to infty, when mu* = infty}
First we study the behavior of $a$, $c_-$ at $\infty$.

\medskip

\proclaim Proposition \label[Proposition:3.1] {\rm (Behavior as $\mu \to \infty$, when $\mu_* = \infty$)}.  
Assume (g0)--(g3) and $\mu_* = \infty$. Then we have
\item{(i)} If $\limsup_{s \to \infty} {G(s)\over s^{2 + {4\over N}}} \leq 0$, then
as $\mu \to \infty$, 
\itemitem{(1)} ${a(\mu)\over \mu} \to \infty$;
\itemitem{(2)} when $N \ge 3$, $c_{-}(\mu) \to \infty$;
\itemitem{(3)} when $N = 2$, assume moreover
$\limsup_{s \to \infty} {g(s)\over s^{1 + {4\over N}}} \le 0$,
then
$c_{-}(\mu) \to \infty$.
\item{(ii)} If
$\liminf_{s \to \infty} {G(s)\over s^{2 + {4\over N}}} = \infty$,
then
    $\displaystyle  {a(\mu)\over \mu} \to 0 \quad \hbox{as}\ \mu \to \infty
    $.


\claim Proof of Proposition \ref[Proposition:3.1].
We use a comparison argument.

\smallskip

\noindent
{\sl 
Step 1: A comparison functional $M(\delta;\mu,u)$.
}

\smallskip

\noindent
For $\delta>0$, we introduce a functional 
$M(\delta;\mu,u) \in C^1((0,\infty)\times E,\R)$ by
    $$  M(\delta;\mu,u) = \half  \norm{\nabla u}_2^2 
        + {\mu\over 2}\norm{u}_2^2
        - {\delta\over 2+{4\over N}} \norm{u}_{2+{4\over N}}^{2+{4\over N}},
    $$
corresponding to the following $L^2$-critical problem:
    $$  -\Delta u + \mu u = \delta u^{1+{4\over N}} \quad \hbox{in}\  \R^N. \eqno\label[3.1]
    $$
Problem \ref[3.1] enjoys the following scaling property.  
Let $\omega_1(x)$ be the least energy solution of \ref[3.1]
with $\mu=1$, $\delta=1$. 
Then
    $$  \omega_{\delta,\mu}(x) = \Bigl({\mu\over \delta}\Bigr)^{N/4} \, \omega_1(\mu^{1/2}x)
    $$
is the least energy solution of \ref[3.1] and its energy level $\beta(\delta;\mu)$
is given by
    $$  \beta(\delta;\mu) = M(\delta;\mu,\omega_{\delta,\mu})
        = \beta_0 \delta^{-N/2} \mu,
    $$
where $\beta_0 = \beta(1;1) > 0$.  We note that $\beta(\delta;\mu)$ is characterized by MP as
    $$  \beta(\delta;\mu) = \inf_{\gamma \in \Gamma_{\delta,\mu}} \max_{\tau \in [0,1]} 
        M(\delta;\mu,\gamma(\tau)),
    $$
where $\Gamma_{\delta,\mu}$ is defined as in Section \ref[Subsection:2.2].

First we consider (i).  

\smallskip

\noindent
{\sl 
Step 2: Under the condition $\limsup_{s \to \infty} {G(s)\over s^{2+{4\over N}}} \le 0$,
${a(\mu)\over \mu} \to\infty$ as $\mu\to\infty$.
}

\smallskip

\noindent
Under the condition in Step 2, for any $\delta > 0$ there exists $A_\delta > 0$ such that
    $$  G(s) \le \half  A_\delta s^2 + {\delta\over 2+{4\over N}} s^{2+{4\over N}}
        \quad \hbox{for}\  s \ge 0. \eqno\label[3.2]
    $$
Thus
    $$  I(\mu,u) \ge M(\delta;\mu-A_\delta,u) \quad
        \hbox{for all}\  (\mu,u) \in (0,\infty)\times E,
    $$
which implies
    $$  a(\mu) \ge \beta(\delta;\mu-A_\delta)
        = \beta_0 \delta^{-N/2} (\mu-A_\delta).
    $$
Thus
    $$  \lim_{\mu \to \infty} {a(\mu)\over \mu} \ge \beta_0 \delta^{-N/2}.
    $$
Since $\delta > 0$ is arbitrary, we have
    $   {a(\mu)\over \mu} \to \infty
    $ as $\mu\to\infty$.

\smallskip

\noindent
{\sl 
Step 3: $c_-(\mu)\to\infty$ as $\mu\to\infty$.
}

\smallskip

\noindent
Next we show $c_-(\mu) \to \infty$ as $\mu \to \infty$.  
Let $u_\mu \in \calC(\mu)$ be the least energy solution with
    $   \half  \norm{u_\mu}_2^2 = c_-(\mu)
    $.
By the Pohozaev identity and the equation, we have
    $$  \eqalignno{
        &{N-2\over 2N} \norm{\nabla u_\mu}_2^2 
            + {\mu\over 2} \norm{u_\mu}_2^2
            = \int_{\R^N} G(u_\mu)\,dx, &\label[3.3]\cr
        &\norm{\nabla u_\mu}_2^2 + \mu \norm{u_\mu}_2^2
            = \int_{\R^N} g(u_\mu) u_\mu \, dx. &\label[3.4]\cr}
    $$
We use \ref[3.3] for $N \ge 3$ and \ref[3.4] for $N=2$.  
We mainly discuss the case $N \ge 3$.  
By \ref[3.2],
    $$  \eqalign{
        {N-2\over 2N} \norm{\nabla u_\mu}_2^2 
            + {\mu\over 2}\norm{u_\mu}_2^2 
        &\le \half  A_\delta \norm{u_\mu}_2^2
            + {\delta\over 2+{4\over N}} \norm{u_\mu}_{2+{4\over N}}^{2+{4\over N}} \cr
        &\le \half  A_\delta \norm{u_\mu}_2^2
            + {\delta C_N\over 2+{4\over N}} \norm{\nabla u_\mu}_2^2 \, \norm{u_\mu}_2^{{4\over N}}. \cr}
    $$
Here we use the Gagliardo-Nirenberg inequality:
    $   \norm{u_\mu}_{2+{4\over N}}^{2+{4\over N}} 
        \le C_N \norm{\nabla u_\mu}_2^2 \, \norm{u_\mu}_2^{{4\over N}}
    $.

We argue indirectly and assume that $c_-(\mu)$ stays bounded for some sequence
$\mu_j\to \infty$.  For the sake of simplicity, we denote the sequence
by \lq\lq$\mu$''.
Thus there exists some constant $B>0$ such that
    $$  c_-(\mu)=\half\norm{u_\mu}_2^2 \le B.
    $$
Choosing $\delta>0$ small so that
    $   {\delta C_N\over 2+{4\over N}} (2B)^{2/N} < {N-2\over 2N}
    $,
we find $\norm{\nabla u_\mu}_2^2$ stays bounded.  

On the other hand, by Step 2, we have
${1\over N} \norm{\nabla u_\mu}_2^2 = a(\mu) \to \infty$, 
which is a contradiction. Hence we have
    $$  c_-(\mu) = \half  \norm{u_\mu}_2^2 \to \infty \quad \hbox{as}\  \mu \to \infty.
    $$
For $N=2$, we use \ref[3.4] instead of \ref[3.3].  
Under the slightly stronger condition\m
    $   \limsup_{s \to \infty} {g(s)\over s^{1+{4\over N}}} \le 0
    $,
we can show $c_-(\mu) \to \infty$ in a similar way.

\smallskip

\noindent
{\sl 
Step 4: Proof of (ii).
}

\smallskip

\noindent
For (ii), under the condition
    $   \liminf_{s \to \infty} {G(s)\over s^{2+{4\over N}}} = \infty
    $,
for any $L>1$ there exists $A_L>0$ such that
    $$  G(s) \ge -\half  A_L s^2 + {L\over 2+{4\over N}} s^{2+{4\over N}} \quad 
        \hbox{for}\ s \ge 0.
    $$
In a similar way, we can show (ii).
\QED


\medskip

\BSScap{\label[Subsection:3.2]. Behavior as $\mu \to \mustarm$, when $\mu_* < \infty$}{\ref[Subsection:3.2]. Behavior as mu to mu*, when mu* < infty}
Next we consider the case $\mu_*<\infty$.

\medskip

\proclaim Proposition \label[Proposition:3.2] 
{\rm (Behavior as $\mu \to \mustarm$, when $\mu_* < \infty$)}.
Assume (g0)--(g3) and $\mu_* < \infty$. Then, as $\mu \to \mustarm$
\item{(1)} $a(\mu) \to \infty$;
\item{(2)} when $N \ge 3$,
$c_{-}(\mu) \to \infty$;
\item{(3)} when $N = 2$, assume moreover
$\limsup_{s \to \infty} {g(s)\over s^{1 + {4\over  N}}} \le 0$.
Then $c_{-}(\mu) \to \infty$.

\medskip


\claim Proof of Proposition \ref[Proposition:3.2].
For $\mu \in (0,\mu_*)$ let $u_\mu \in \calC(\mu)$ be a least energy
solution with $\half  \norm{u_\mu}_2^2 = c_-(\mu)$.
We note that $a(\mu) = I(\mu,u_\mu) = {1\over N} \norm{\nabla u_\mu}_2^2$.

Arguing indirectly, we assume that
$a(\mu) = {1\over N} \norm{\nabla u_\mu}_2^2$ stays bounded along a sequence
$\mu_j\to \mustarm$.  That is, for some $B>0$, 
    $$  \norm{\nabla u_{\mu_j}}_2^2 \le B.
    $$
For notational simplicity, we simply denote the sequence as $\mu$.
Proof consists of several steps.

\smallskip

\noindent
{\sl
Step 1. $\norm{u_\mu}_2^2$ stays bounded as $\mu \to \mustarm$.
}

\smallskip

\noindent
We argue indirectly and assume $\norm{u_{\mu_j}}_2^2 \to \infty$
for some sequence $\mu_j \to \mustarm$.
Again we simply denote the sequence as $\mu$.
We set
    $$  v_\mu(x) = u_\mu(x/t_\mu), \ \hbox{where} \ 
        t_\mu = \norm{u_\mu}_2^{-2/N} \to 0 \quad \hbox{as}\  \mu \to \mustarm.
    $$
We have for $\mu$ close to $\mu_*$
    $$  \eqalignno{
        &\norm{\nabla v_\mu}_2^2 = t_\mu^{N-2} \norm{\nabla u_\mu}_2^2
            \le \norm{\nabla u_\mu}_2^2 \le B, \cr
        &\norm{v_\mu}_2^2 = t_\mu^N \norm{u_\mu}_2^2 = 1. &\label[3.5]\cr}
    $$
Thus $(v_\mu)$ stays bounded in $E$ and we may assume
$v_\mu \wlimit v_0$ weakly in $E$ as $\mu \to \mustarm$.

Since $v_\mu$ solves
    $   -t_\mu^2 \Delta v_\mu + \mu v_\mu = g(v_\mu)$ in $\R^N$,
the weak limit $v_0 \in E$ satisfies
    $$  \mu_* v_0 = g(v_0) \quad \hbox{in}\  \R^N.
    $$
By the behavior $g(s) = o(s)$ near $0$, we conclude $v_0=0$.
Compactness of the embedding $E \subset L^r(\R^N)$
for $r \in (2,2^*)$ implies $v_\mu \to 0$ strongly in $L^r(\R^N)$
for $r \in (2,2^*)$.
Thus
    $$  t_\mu^2\norm{\nabla v_\mu}_2^2 + \mu \norm{v_\mu}_2^2
        = \int_{\R^N} g(v_\mu)v_\mu \, dx \to 0 \quad \hbox{as}\ \mu \to \mustarm.
    $$
In particular, $\norm{v_\mu}_2^2 \to 0$, which is
in contradiction to \ref[3.5]. Thus $(u_\mu)$ is bounded in $L^2(\R^N)$
and thus $(u_\mu)$ is bounded in $E$.
We may assume $u_\mu \wlimit u_0$ weakly in $E$
as $\mu \to \mustarm$.

\smallskip

\noindent
{\sl Step 2: $u_0=0$ and $u_\mu \to 0$ strongly in $E$.
}

\smallskip 

\noindent
Since $u_\mu$ satisfies
    $-\Delta u_\mu + \mu u_\mu = g(u_\mu)$ in $\R^N$,
the limit function $u_0 \in E$ solves
    $$  -\Delta u_0 + \mu_* u_0 = g(u_0) \quad \hbox{in}\  \R^N.
    $$
Thus by the Pohozaev identity,
    $$  {N-2\over 2N} \norm{\nabla u_0}_2^2 
            + \int_{\R^N}{\mu_*\over 2} u_0^2 -  G(u_0)\, dx = 0.
        \eqno\label[3.6]
    $$
Since
    $$  \eqalign{
        &{\mu_*\over 2}s^2 - G(s) > 0 \quad \hbox{for}\  s \in (0,\delta), \cr
        &{\mu_*\over 2}s^2 - G(s) \ge 0 \quad \hbox{for all}\  s \ge 0\cr}
    $$
for some $\delta>0$, \ref[3.6] implies $u_0=0$.

Thus $u_\mu \wlimit 0$ weakly in $E$ and strongly
in $L^r(\R^N)$ for $r \in (2,2^*)$, which implies
    $$  \norm{\nabla u_\mu}_2^2 + \mu \norm{u_\mu}_2^2
        = \int_{\R^N} g(u_\mu)u_\mu\, dx \to 0.
    $$
Thus $u_\mu \to 0$ strongly in $E$.

\smallskip

\noindent
{\sl 
Step 3: $a(\mu) \to \infty$ as $\mu \to \mustarm$.
}

\smallskip

\noindent
By the monotonicity of $\mu\mapsto a(\mu)$, we have
    $$  0<a\left({\mu_*\over 2}\right) \leq a(\mu) ={1\over N}\norm{\nabla u_\mu}_2^2 
        \quad \hbox{for}\ \mu\in \left[{\mu_*\over 2}, \mu_*\right),
    $$
which is incompatible with Step 2.
Thus we have Step 3.

\smallskip

\noindent
{\sl 
Step 4: $c_-(\mu) \to \infty$ as $\mu \to \mustarm$.
}

\smallskip

\noindent
First we consider the case $N \ge 3$.  
By the Pohozaev identity we have
    $$  {N-2\over 2} a(\mu)
        = {N-2\over 2N} \norm{\nabla u_\mu}_2^2
        = -{\mu\over 2}\norm{u_\mu}_2^2 + \int_{\R^N} G(u_\mu)\,dx.
    $$
Thus by \ref[2.2]
    $$  {N-2\over 2} a(\mu) \le \half  (\mu_* - \mu)\norm{u_\mu}_2^2
        = (\mu_* - \mu) c_-(\mu).
    $$
Hence by Step 3 we have $c_-(\mu) \to \infty$ as $\mu \to \mustarm$.

Next we consider the case $N=2$.  
By the assumption
    $   \limsup_{s \to \infty} {g(s)\over s^{1+4/N}} \le 0
    $,
for any $\epsilon > 0$ there exists $A_\epsilon > 0$ such that
    $$  g(s)s \le \epsilon s^4 + A_\epsilon s^2, \quad s \ge 0.
    $$
Thus
    $$  \eqalign{
        \norm{\nabla u_\mu}_2^2 + \mu \norm{u_\mu}_2^2
        &= \int_{\R^2} g(u_\mu) u_\mu \, dx 
        \le \epsilon \norm{u_\mu}_4^4 + A_\epsilon \norm{u_\mu}_2^2 \cr
        & \le \epsilon C_2\norm{\nabla u_\mu}_2^2 \norm{u_\mu}_2^2
            + A_\epsilon \norm{u_\mu}_2^2. \cr}
    $$
Here we use Gagliardo-Nirenberg inequality:
$\norm{u}_4^4 \le C_2 \norm{\nabla u}_2^2 \norm{u}_2^2$.

If $c_-(\mu)$ stays bounded as $\mu \to \mustarm$, i.e.,
$\norm{u_\mu}_2^2 \le M$ for some constant $M>0$ independent of $\mu$,
we have
    $$  \norm{\nabla u_\mu}_2^2 
        \le \epsilon C_2 \norm{\nabla u_\mu}_2^2 M + A_\epsilon M.
    $$
Choosing $\epsilon>0$ such that $\epsilon C_2 M < 1$, we observe that
    $   a(\mu) = {1\over N} \norm{\nabla u_\mu}_2^2
    $
stays bounded as $\mu \to \mustarm$, which contradicts Step 3.
Thus
    $   c_-(\mu) \to \infty \quad \hbox{as}\  \mu \to \mustarm
    $.
\QED


\medskip

\BSScap{\label[Subsection:3.3]. Behavior as $\mu \to 0^+$}{\ref[Subsection:3.3]. Behavior as mu to 0^+}
Finally we study the behavior as $\mu \to 0^+$.  
Here the behavior of $g(s)$ as $s\sim 0$ is important.
We assume

\smallskip

{\parindent=1.5\parindent

\item{(g1')} For some $p_0>1$, $\lim_{s\to 0^+} {g(s)\over s^{p_0}} =1$.

}

\smallskip

\noindent
We note (g1) and (g3) follow from (g1').

\medskip

\proclaim Proposition \label[Proposition:3.3]. (Behavior as $\mu \to 0^+$).
Assume (g0), (g1') and (g2). Then we have
\item{(i)} if $p_0$ in (g1') satisfies $p_0 \in (1, 2^*-1)$, then as $\mu \to 0^+$
    $$  a(\mu) = \mu^{\beta+1}\,(a_{0,{p_0}} + o(1)), \qquad
        c_\pm(\mu) = \mu^\beta (c_{0,{p_0}}+o(1)), 
    $$
where
    $$  \eqalignno{
        \beta &= {2\over p_0 - 1} - {N\over 2} \in (-1,+\infty), &\label[3.7]\cr
        a_{0,{p_0}} &= \half  \norm{\nabla \omega_{p_0}}_2^2 + \half  \norm{\omega_{p_0}}_2^2
            - {1\over {p_0}+1}\norm{\omega_{p_0}}_{{p_0}+1}^{{p_0}+1} \in (0,\infty), \cr
        c_{0,{p_0}} &= \half  \norm{\omega_{p_0}}_2^2. \cr}
    $$
Here $\omega_{p_0}(x)$ is the unique ground state solution of
    $   -\Delta \omega + \omega = \omega^{p_0} \quad \hbox{in}\  \R^N
    $.
\item{(ii)} When $N \ge 3$, if $p_0 \in [2^*-1,+\infty)$, then
    $$  \lim_{\mu \to 0^+} a(\mu) \in (0,\infty).
    $$

\medskip

\claim Proof of Proposition \ref[Proposition:3.3]. \m
First we consider (i) and suppose $p_0 \in (1,2^*-1)$.
Let $u_\mu(x) \in \calC(\mu)$ be the least energy solution.
We introduce $w_\mu(x)$ by
    $$  u_\mu(x) = \mu^{1/(p_0-1)} w_\mu(\mu^{1/2}x).
    $$
We also introduce $h(s) = g(s) - s^{p_0}$ and $H(s) = \int_0^s h(t)\,dt$.
We observe that $h(s)$ and $H(s)$ satisfy 
$h(s) = o(s^{p_0})$, $H(s) = o(s^{p_0+1})$
as $s \to 0^+$.

We compute
    $$  I(\mu,u_\mu) = \half  \norm{\nabla u_\mu}_2^2 + {\mu\over 2} \norm{u_\mu}_2^2
            - {1\over p_0+1} \norm{u_\mu}_{p_0+1}^{p_0+1}
            - \int_{\R^N} H(u_\mu)\,dx 
        = \mu^{\beta+1} L(\mu,w_\mu), 
    $$
where
    $$  L(\mu,w) = \half  \norm{\nabla w}_2^2
            + \half  \norm{w}_2^2 - {1\over p_0+1}\norm{w}_{p_0+1}^{p_0+1}
            - \mu^{-{p_0+1\over p_0-1}} \int_{\R^N} H(\mu^{1/(p_0-1)}w)\,dx.
    $$
It is easily seen that $w_\mu$ is a MP critical point of $w\mapsto L(\mu,w)$
and $w_\mu \to \omega_{p_0}$ strongly in $E$ as $\mu \to 0$,
where $\omega_{p_0}$ is the least energy solution of
    $   -\Delta \omega + \omega = \omega^{p_0}$ in $\R^N$.
Thus we have
    $$  \eqalign{
        a(\mu) &= {1\over N} \norm{\nabla u_\mu}_2^2
            = {1\over N} \mu^{\beta+1} \norm{\nabla w_\mu}_2^2
            = \mu^{\beta+1} ({1\over N} \norm{\nabla \omega_{p_0}}_2^2 + o(1)), \cr
        \half  \norm{u_\mu}_2^2
        &= \half  \mu^\beta \norm{w_\mu}_2^2
        = \mu^\beta (\half\norm{\omega_{p_0}}_2^2 + o(1)).\cr}
    $$
Since $u_\mu\in \calC(\mu)$ is arbitrary,  we have (i).

\smallskip

Next we prove (ii).  Suppose $N \ge 3$ and $p_0 \ge 2^*-1$.
By assumption (g2), there exists $C>0$ such that
    $$  G(s) \le C s^{2^*}  \quad \hbox{for}\ s \ge 0.
    $$
Thus by Sobolev embedding $\calD^{1,2}(\R^N) \subset L^{2^*}(\R^N)$,
    $$  I(\mu,u) \ge \half  \norm{\nabla u}_2^2 + {\mu\over 2} \norm{u}_2^2
            - C \norm{u}_{2^*}^{2^*} 
        \ge \half  \norm{\nabla u}_2^2 - C'\norm{\nabla u}_2^{2^*}. 
    $$
We observe that $L_*(u) = \half  \norm{\nabla u}_2^2 - C'\norm{\nabla u}_2^{2^*}$
has a MP geometry in $\calD^{1,2}(\R^N)$ and 
    $$  a(\mu) \ge (\hbox{MP minimax value for $L_*$}) >0 
        \quad \hbox{for all} \ \mu\in (0,\mu_*).
    $$
Thus (ii) holds.
\QED


\medskip

\noindent
We note that for $\beta$ given in \ref[3.7] satisfies
$\beta>0$ for $p_0 \in (1,1+{4\over N})$ and 
$\beta \in (-1,0)$ for $p_0 \in (1+{4\over N}, 2^*-1)$.
Thus we have

\medskip

\proclaim Corollary \label[Corollary:3.4].
Assume (g$_0$), (g1'), (g2) and (g3).  Then
\item{(i)} if $p_0 \in (1,1+{4\over N})$, then
    $a(\mu) \to 0$, ${a(\mu)\over \mu} \to 0$, $c_\pm(\mu) \to 0$ as $\mu\to 0^+$;
\item{(ii)} If $p_0 \in (1+{4\over N}, 2^*-1)$, then
    $a(\mu) \to 0$,  ${a(\mu)\over \mu} \to \infty$, $c_\pm(\mu) \to \infty$
as $\mu \to 0^+$.
\QED

\medskip

\BScap{\label[Section:4]. $(PSPC)$ condition for $\widetilde{J}_m(\lambda,u)$}{\ref[Section:4]. $(PSPC)$ condition for widetilde J_m}

\BSS{\label[Subsection:4.1].  Proposition \ref[Proposition:4.1]}
To apply Theorem \ref[Theorem:1.6], we need to verify the $(PSPC)$ condition.  
Here we consider $(PSPC)$ under the following condition:

{\parindent=1.5\parindent
\smallskip

{\item{(p1)} $p_0$ in (g1') satisfies 
    $\displaystyle 1 < p_0 <2^*-1$.  When $N\geq 3$, if $p_0>{N\over N-2}$, we assume
    in addition
        $$  g(s)s \leq 2^* G(s) \quad \hbox{for}\ s\geq 0.      \eqno\label[4.1]
        $$
\smallskip
}}

\noindent
In other words, we assume

\bitem when $N=2,3$, we assume (p0) in Theorems \ref[Theorem:1.1] and \ref[Theorem:1.2], that is, 
$1+{4\over N}<p_0\leq {N\over N-2}$

\noindent or

\bitem when $N\geq 3$, for ${N\over N-2}<p_0<2^*-1$ we assume the global condition
\ref[4.1].

\smallskip

\noindent
We note that \ref[4.1] is also found in [\cite[JZZ:31], \cite[MS:33]]
and it implies that $s\mapsto {G(s)\over s^{2^*}}$ is non-increasing in $(0,\infty)$.

We have the following.

\medskip

\proclaim Proposition \label[Proposition:4.1].  
Suppose (g0), (g1'), (g2) and (p1).
Moreover assume one of the following three conditions:
{\parindent=1.5\parindent
\item{(b1)} $\mu_* < \infty$;
\item{(b2)} $\mu_* = \infty$ and there exist $\theta \in (2, 2 + {4\over N})$ and 
$L \ge 1$ such that
    $$  \theta G(s) \ge g(s)s \quad \hbox{for } s \ge L;
    $$
\item{(b3)} $\mu_* = \infty$ and there exist $\theta \in (2 + {4\over N}, 2^*)$ and 
$L \ge 1$ such that
    $$  0 < \theta G(s) \le g(s)s \quad \hbox{for } s \ge L.
    $$
}
Then for $m > 0$, $\widetilde{J}_m(\lambda,u)$ satisfies the $(PSPC)_b$ condition 
for any $b \ne 0$.

\medskip


To prove Proposition \ref[Proposition:4.1], 
for $m > 0$, $b\in\R$ we suppose that 
$(\lambda_j,u_j)_{j=1}^\infty \subset \R \times E$ is a $(PSPC)_b$ sequence
and we consider three cases separately;  

\smallskip

{\parindent=2.5\parindent
\item{\bf Case A:} $\lambda_j \to \lambda_0$ as $j \to \infty$ 
for some $\lambda_0 \in \R$;
\item{\bf Case B:} $\lambda_j \to \infty$ as $j \to \infty$;
\item{\bf Case C:} $\lambda_j \to -\infty$ as $j \to \infty$. 
}

\smallskip

\noindent
In what follows, we use notation:
    $$  \mu_j = e^{\lambda_j}.
    $$
Case A is addressed in Section \ref[Subsection:2.2] a) essentially.  

\medskip

\proclaim Lemma \label[Lemma:4.2].
Suppose that Case A occurs.  Then $(\lambda_j,u_j)_{j=1}^\infty$ has 
a strongly convergent subsequence in $\RE$.

\claim Proof (Proof of Proposition \ref[Proposition:4.1]: Case A).  
For a $(PSPC)_b$ sequence $(\lambda_j,u_j)_{j=1}^\infty$, we may assume 
$\mu_j\to\mu_0\in (0,\infty)$.  Then $(\mu_j,u_j)_{j=1}^\infty$, $\mu_j=e^{\lambda_j}$
satisfies $I(\mu_j,u_j)\to b+\mu_0 m$, $\norm{\partial_u I(\mu_j,u_j)}_{E^*}\m\to 0$,
$P(\mu_j,u_j)\to 0$.  Thus Lemma \ref[Lemma:4.2] holds by Proposition \ref[Proposition:2.7].
\QED

\medskip

To deal with Cases B and C, first we note that $(\lambda_j,u_j)_{j=1}^\infty$ satisfies
\ref[2.11]--\ref[2.14],
from which we have
    $$  \eqalignno{
        &\mu_j \norm{u_j}_{2}^{2} = 2 \mu_j m + o(1), &\label[4.2]\cr
        &\norm{\nabla u_j}_{2}^{2} = N m \mu_j + N b + o(1), &\label[4.3]\cr
        &\int_{\R^N} G(u_j)\, dx = {N\over 2} m \mu_j + {N-2\over 2} b + o(1), 
            &\label[4.4]\cr
        &\int_{\R^N} g(u_j)u_j\, dx = (N+2)m\mu_j + N b + o(1). &\label[4.5]\cr}
    $$
We also note 
the following properties of $(PSPC)_b$ sequence.

\medskip

\proclaim Lemma \label[Lemma:4.3].  
For a $(PSP)_b$ sequence $(\lambda_j,u_j)_{j=1}^\infty$, we have
\item{(i)} $\norm{u_{j-}}_E \to 0$, $(1+\mu_j) \norm{u_{j-}}_2^2 \to 0$.
\flushright{\label[4.6]}
\item{(ii)} $(\lambda_j,u_{j+})_{j=1}^\infty$ satisfies
    $$  \eqalignno{
        &\norm{u_j}_E^2 = \norm{u_{j+}}_E^2+o(1), &\label[4.7]\cr
        &\norm{\nabla u_j}_2^2 = \norm{\nabla u_{j+}}_2^2 + o(1),
            &\label[4.8]\cr
        &\mu_j \norm{u_j}_2^2 = \mu_j \norm{u_{j+}}_2^2 + o(1),   &\label[4.9]\cr
        &\int_{\R^N} G(u_j)\,dx =\int_{\R^N} G(u_{j+})\,dx + o(1), &\label[4.10]\cr
        &\int_{\R^N} g(u_j)u_j\,dx = \int_{\R^N} g(u_{j+})u_{j+}\,dx + o(1). &\label[4.11]\cr}
    $$

\medskip

\claim Proof.
(i) For a $(PSPC)_b$ sequence $(\lambda_j,u_j)_{j=1}^\infty$, we have\m
    $   \abs{\partial_u \widetilde{J}_m(\lambda_j,u_j) u_{j-}} 
        \leq \norm{\partial_u \widetilde{J}_m(\lambda_j,u_j)}_{E^*} \norm{u_j}_E 
        \to 0
    $.
By \ref[1.3],  
    $$  \norm{\nabla u_{j-}}_2^2 + (1+\mu_j)\norm{u_{j-}}_2^2 
        = -\partial_u\widetilde J_m(\lambda_j,u_j)u_{j-} \to 0.
    $$  
Thus (i) holds.
For (ii) we observe that \ref[4.7]--\ref[4.9] follow from (i).
It is also easy to see that
    $$  \eqalign{
        &\int_{\R^N} g(u_j)u_j\,dx = \int_{\R^N} g(u_{j+})u_{j+}\,dx 
            + \int_{\R^N} g(-u_{j-})(-u_{j-})\,dx \cr
        &\qquad =\int_{\R^N} g(u_{j+})u_{j+}\,dx + o(1) \cr}
    $$
and thus \ref[4.11] follows.  \ref[4.10] can be shown in a similar
way.  Here we use the conditions (g1)--(g2).  \QED


\medskip

\BSS{\label[Subsection:4.2]. Proof of Proposition \ref[Proposition:4.1]: Case B}
Next we consider Case B.

\medskip

\proclaim Lemma \label[Lemma:4.4].  
Assume (b1), i.e., $\mu_* < \infty$. Then Case B cannot take place.

\medskip


\claim Proof.
Suppose that $(\lambda_j, u_j)_{j=1}^\infty$ is a $(PSP)_b$ sequence 
with $\mu_j \to \infty$.
Under the assumption $\mu_* < \infty$, we have \ref[2.2]. Thus by
Lemma \ref[Lemma:4.3], \ref[4.4],
    $$  \eqalign{
        {\mu_*\over 2}\norm{u_j}_2^2
        &= {\mu_*\over 2}\norm{u_{j+}}_2^2 + o(1) 
        \geq \int_{\R^N} G(u_{j+})\,dx + o(1) \cr
        &= \int_{\R^N} G(u_j)\,dx + o(1) 
        ={N\over 2}\mu_j m + {N-2\over 2} b + o(1). \cr}
    $$
Since $\mu_j \to \infty$, $\norm{u_j}_2^2\to\infty$.
On the other hand, \ref[4.2] implies 
    $\norm{u_j}_2^2 \to 2m$.
This is a contradiction and 
$\mu_j \to \infty$ cannot take place.
\QED

\medskip

\proclaim Lemma \label[Lemma:4.5].  
Assume (b2) or (b3) in Proposition \ref[Proposition:4.1]. Then Case B cannot take place.

\medskip


\claim Proof.
Suppose that $\mu_j \to \infty$ occurs for a $(PSPC)_b$ sequence 
$(\lambda_j, u_j)_{j=1}^\infty$.  
The proof consists of two steps.

\smallskip

\noindent
{\sl Step 1. For $L \ge 1$, we set
    $   D(L,j) = \{x \in \R^N \,;\, u_j(x) \ge L\}
    $.
Then for any $L\geq 1$,
    $$  \eqalignno{
        &{1\over \mu_j}\int_{D(L,j)} G(u_j)\,dx \;\to\; {N\over 2}m, &\label[4.12]\cr
        &{1\over \mu_j}\int_{D(L,j)} g(u_j)u_j\,dx \;\to\; (N+2)m. &\label[4.13]\cr}
    $$
}

\noindent
We show \ref[4.12]. \ref[4.13] can be shown in a similar way.
For any $L\geq 1$ there exists a constant $A_L>0$ such that
    $$  \abs{G(s)} \le A_L s^2 \quad \hbox{for}\  s \in (-\infty,L].
    $$
Thus by \ref[4.2]
    $$  {1\over \mu_j}\int_{\R^N \setminus D(L,j)} \abs{G(u_j)}\,dx
        \le {1\over \mu_j} A_L \norm{u_j}_2^2
        = {1\over \mu_j} A_L (2m+o(1)) \to 0.
        \eqno\label[4.14]
    $$
On the other hand, by \ref[4.4] we have
    ${1\over \mu_j}\int_{\R^N} G(u_j)\,dx \;\to\; {N\over 2} m$.
\ref[4.12] follows from \ref[4.14].

\smallskip

\noindent
{\sl Step 2. Under (b2) or (b3), $\mu_j \to \infty$ cannot take a place.
}

\smallskip

\noindent
We consider (b2). We can deal with (b3) in a similar way.
Let $\theta$, $L$ be given in (b2).  We have
    $$  \theta\, G(u_j(x)) \ge g(u_j(x))\,u_j(x) \quad \hbox{for}\  x \in D(L,j),
    $$
from which we have
    $$  \theta \lim_{j\to\infty}{1\over \mu_j}\int_{D(L,j)} G(u_j)\,dx
        \;\ge\; \lim_{j\to\infty}{1\over \mu_j}\int_{D(L,j)} g(u_j)u_j\,dx.
    $$
Thus by Step 1,
    $   \theta {N\over 2}m \ge (N+2)m
    $.
This is incompatible with $\theta \in (2,2+{4\over N})$.  
Therefore $\mu_j \to \infty$ cannot take a place.
\QED

\medskip


By Lemmas \ref[Lemma:4.4] and \ref[Lemma:4.5] under the assumption of Proposition \ref[Proposition:4.1], 
Case B cannot take place.

\medskip

\BSS{\label[Subsection:4.3]. Proof of Proposition \ref[Proposition:4.1]: Case C}
Finally we deal with Case C.
Here condition (p1) plays a role.

\medskip

\proclaim Lemma \label[Lemma:4.6].  
Assume (g0), (g1'), (g2) and (p1).
Then for $b \ne 0$ Case C cannot take place.

\medskip

Lemma \ref[Lemma:4.6] in a special case $p_0 = 1 + {4\over N}$ 
is shown in [\cite[CGIT:15]].
In this section we give a proof of Lemma \ref[Lemma:4.6] using ideas in [\cite[CGIT:15]].

\medskip

\noindent
{\bf a) Function space $F_{p_0+1}$
}

\medskip

\noindent
For $p_0\in (1,2^*-1)$ we define a function space $F_{p_0+1}$ by
    $$  F_{p_0+1} = \{ u \in L^{p_0+1}(\R^{N}) \,;\, \nabla u \in L^{2}(\R^{N}),\,
        u(x) = u(\abs{x}) \}
    $$
and we equip a norm
    $$  \norm{u}_{F_{p_0+1}} = \norm{\nabla u}_{2} + \norm{u}_{p_0+1}.
    $$
$F_{p_0+1}$ is a natural space to study the zero-mass problem
    $$  - \Delta u = g(u) \quad \hbox{in}\  \R^{N}, \eqno\label[4.15]
    $$
which appears as the limit equation of our problem
$- \Delta u + \mu u = g(u)$ as $\mu \to 0$.

We observe $F_{p_0+1} \subset C(\R^N \setminus \{0\},\R)$ and 
$F_{p_0+1} \subset L^q(\R^N)$ for $q \in [p_0+1, 2^*]$ when $N \geq 3$, 
for $q\in [p_0+1,\infty)$ when $N=2$. 
The embeddings $F_{p_0+1} \subset L^q(\R^N)$, 
$F_{p_0+1} \subset L^r_{loc}(\R^N)$ are compact for $q\in (p_0+1, 2^*)$ 
and $r \in [1,2^*)$.  

\medskip

We have the following Liouville type result.

\medskip

\proclaim Proposition \label[Proposition:4.7].  
Let $N \geq 2$ and suppose (g0), (g1') and 
$1<p_0\leq {N\over N-2}$ for $N\geq 3$ and $1<p_0<\infty$ for $N=2$.
Then 
the zero-mass problem \ref[4.15] has no positive solutions in $F_{p_0+1}$. 
\QED

\medskip

Proposition \ref[Proposition:4.7] is essentially due to Armstrong-Sirakov [\cite[AS:1]]
and it is shown in [\cite[CGIT:15], Theorem \raw[2.6]].

\medskip

We also have the following result under \ref[4.1]:

\medskip

\proclaim Proposition \label[Proposition:4.8].
Let $N \ge 3$ and suppose (g0), (g1') with $p_0\in (1,2^*-1)$, (g2) and \ref[4.1].
Then the zero-mass problem \ref[4.15] has no positive
solutions in $F_{p_0+1}$. \QED

\smallskip

Proposition \ref[Proposition:4.8] is shown for $p_0 = 1 + {4\over N}$
in [\cite[CGIT:15], Proposition \raw[1.2] (ii)].
We can follow the argument in [\cite[CGIT:15]].  For a solution $u$ 
of \ref[4.15], by Pohozaev identity and $\norm{\nabla u}_2^2 = \intRN g(u)u\, dx$, 
we have
    $$  {N-2\over 2N}\intRN g(u)u\, dx -\intRN G(u) \, dx =0.
    $$
The conditions \ref[4.1] and (g1') with $p_0\in (1,2^*-1)$ imply
that $2^*G(s)-g(s)s$ is positive for small $s>0$ and non-negative 
for all $s\geq 0$.  
Thus $u(x)=0$ for all $x\in \R^N$.  We also note that for $N=2$, 
\ref[4.15] has no positive
solutions in $F_{p_0+1}$ under the condition $G(s)\geq 0$ for all $s\geq 0$.

As a corollary to Propositions \ref[Proposition:4.7] and \ref[Proposition:4.8],
we have

\medskip

\proclaim Corollary \label[Corollary:4.9].
Let $N \ge 2$ and assume (g0), (g1'), (g2) and (p1).
Then the zero-mass problem \ref[4.15] has no positive
solutions in $F_{p_0+1}$.
\QED

\medskip

We also have the following radial lemma.

\medskip

\proclaim Lemma \label[Lemma:4.10].  
\item{(i)} ([\cite[BL:6], Radial Lemma \raw[A.III]])
When $N \geq 3$, there exists a constant $C_{N} > 0$ such that 
for $u \in \calD^{1,2}_{r}(\R^{N})$
    $$  \abs{u(x)} \leq C_{N} \norm{\nabla u}_{2} 
        \abs{x}^{-{N-2\over 2}} \quad \hbox{for}\  \abs{x} \geq 1.
    $$
\item{(ii)} When $N=2$, there exists a constant $C_{2}>0$ such that 
for $u \in F_{p_0+1}$
    $$  \abs{u(x)}^{{p_0+3\over 2}} 
        \leq C_{2} \norm{\nabla u}_{2} \norm{u}_{p_0+1}^{{p_0+1\over 2}} 
        \abs{x}^{-1} \quad \hbox{for}\ x \neq 0. \eqno\label[4.16]
    $$

\medskip

\claim Proof.  (c.f. Lemma \raw[2.4] in [\cite[CGIT:15]]).
We prove (ii). It suffices to show \ref[4.16] for $u \in C_{0,r}^{\infty}(\R^{2})$.  
We write $r=\abs{x}$ and compute
    $$  \eqalign{
    {d\over dr}\left( r \abs{u(r)}^{{p_0+3\over 2}} \right) 
    &= \abs{u(r)}^{{p_0+3\over 2}} + {p_0+3\over 2} r \abs{u(r)}^{{p_0+1\over 2}} 
        u_{r}(r)\, \hbox{sgn}(u(r)) \cr
    &\geq - {p_0+3\over 2} r \abs{u(r)}^{{p_0+1\over 2}} \abs{u_{r}(r)}.\cr}
    $$
Integrating on $[r,\infty)$, we have
    $$  \eqalign{
    r \abs{u(r)}^{{p_0+3\over 2}}
    &= - \int_{r}^{\infty} {d\over ds}\left( s \abs{u(s)}^{{p_0+3\over 2}} \right) ds
    \leq {p_0+3\over 2} \int_{r}^{\infty} 
        s \abs{u(s)}^{{p_0+1\over 2}} \abs{u_{r}(s)} ds \cr
    &\leq {p_0+3\over 2} \norm{u}_{p_0+1}^{{p_0+1\over 2}} \norm{\nabla u}_{2}. \cr}
    $$
Thus we get \ref[4.16]. 
\QED

\medskip

\noindent
{\bf b) Boundedness of $(PSPC)_b$ sequence in $F_{p_0+1}$
}

\medskip

\noindent
We show the boundedness of a $(PSPC)_b$ sequence with $\mu_{j} \to 0$.
By \ref[4.3], we may assume $b\geq 0$.

\medskip

\claim Lemma \label[Lemma:4.11].  
Assume (g0), (g1'), (g2) with $p_0\in (1,2^*-1)$.
Suppose that $(\lambda_{j}, u_{j})_{j=1}^{\infty}$ is a 
$(PSPC)_b$ sequence with $b\not=0$ and $\mu_{j}\to 0$.  
Then $(u_{j})_{j=1}^{\infty}$ is bounded in $F_{p_0+1}$.

\medskip

\claim Proof.  
We may assume $b>0$.  
By \ref[4.3], we observe that $\norm{\nabla u_{j}}_{2}^{2}$ is bounded.  
Thus we need to show boundedness of $\norm{u_{j}}_{p_0+1}$.  
We argue indirectly and assume
    $   t_{j} = \norm{u_{j}}_{p_0+1}^{-{p_0+1\over N}} \to 0
    $.
We set
    $$  v_{j}(x) = u_{j}(x/t_{j}).
    $$
We observe from \ref[4.2]--\ref[4.5]
    $$  \eqalignno{
    &\norm{\nabla v_{j}}_{2}^{2} = t_{j}^{N-2} \norm{\nabla u_{j}}_{2}^{2}
        \to \cases{ Nb  & when $N=2$, \cr
                    0   & when $N \geq 3$,\cr}
    &\label[4.17]\cr
    &\norm{v_{j}}_{p_0+1}^{p_0+1} 
        = t_{j}^{N} \norm{u_{j}}_{p_0+1}^{p_0+1} = 1, &\label[4.18]\cr
    &\mu_{j}\norm{v_{j}}_{2}^{2} = t_{j}^{N} \mu_{j} \norm{u_{j}}_{2}^{2} 
        \to 0, &\label[4.19]\cr
    &\int_{\R^{N}} g(v_{j}) v_{j} \,dx 
        = t_{j}^{N} \int_{\R^{N}} g(u_{j}) u_{j} \,dx \to 0. &\label[4.20]\cr}
    $$
We also observe from \ref[4.6] that
    $$  \norm{v_{j-}}_E^2 = \norm{\nabla v_{j-}}_2^2 + \norm{v_{j-}}_2^2
        = t_j^{N-2} \norm{\nabla u_{j-}}_2^2 + t_j^N \norm{u_{j-}}_2^2 \to 0.
        \eqno\label[4.21]
    $$
The proof consists of two steps.

\smallskip

\noindent
{\sl
Step 1. $v_{j} \wlimit 0$ weakly in $F_{p_0+1}$.
}

\smallskip

\noindent
By \ref[4.17], \ref[4.18], $(v_{j})_{j=1}^{\infty}$ is bounded 
in $F_{p_0+1}$ and, after extracting a subsequence,  
we may assume $v_{j} \wlimit v_{0}$ weakly 
in $F_{p_0+1}$ for some $v_{0} \in F_{p_0+1}$.

\smallskip

\noindent
When $N \geq 3$, from \ref[4.17] we observe $\nabla v_{0}=0$.  
Thus $v_{0}=0$ follows.  
We show $v_{0}=0$ also for $N=2$.

Let $N=2$, by the Gagliardo-Nirenberg inequality we have
    $$  t_{j}^{-2} = \norm{u_{j}}_{p_0+1}^{p_0+1} 
        \leq C \norm{\nabla u_{j}}_{2}^{p_0-1\over 2} \norm{u_{j}}_{2}^{2}
        = C (Nb + o(1))^{p_0-1\over 2} \norm{u_{j}}_{2}^{2}.
    $$
Thus for some constant $c>0$
    $$ t_{j}^{-2} \leq c \norm{u_{j}}_{2}^{2} \leq c \norm{u_{j}}_E^{2}.
    $$
For any given $\varphi \in E$, we set $\psi_j(x)=\varphi(t_{j}x)$.  We have
for some constant $C_{\varphi}>0$ depending on $\varphi$
    $$  \norm{\psi_j}_{E}
        = (\norm{\nabla \varphi}_{2}^{2} + t_{j}^{-2}\norm{\varphi}_{2}^{2})^{1/2}
        \leq \norm\varphi_E ( 1+t_j^{-1})   \leq C_\varphi(1+\norm{u_j}_E). 
    $$
Since $(\lambda_{j}, u_{j})_{j=1}^{\infty}$ is a $(PSPC)_b$ sequence,
    $$  \eqalign{
    \abs{ \partial_{u}\widetilde{J}_{m}(\lambda_{j},u_{j}) \psi_j} 
    &\leq 
        \norm{\partial_{u}\widetilde{J}_{m}(\lambda_{j},u_{j})}_{E^{*}}
        \norm{\psi_j}_{E} \cr
    &\leq C_{\varphi}(1+\norm{u_{j}}_{E})
        \norm{\partial_{u}\widetilde{J}_{m}(\lambda_{j},u_{j})}_{E^{*}}
        \to 0. \cr}
    $$
Thus
    $$  \eqalign{
    &t_{j}^{2} (\nabla v_{j}, \nabla \varphi)_{2} + \mu_{j}(v_{j}, \varphi)_{2} 
        - \int_{\R^{N}} g(v_{j})\varphi \,dx\cr
    =& t_{j}^{2}\left( (\nabla u_{j}, \nabla \psi_j)_{2} 
        + \mu_{j}(u_{j}, \psi_j)_{2} - \int_{\R^{N}} g(u_{j})\psi_j \,dx\right)\cr
    =& t_{j}^{2} \partial_{u}
        \widetilde{J}_{m}(\lambda_{j}, u_{j}) \psi_j \to 0, \cr}
    $$
from which we have
    $$  \int_{\R^{N}} g(v_{0}) \varphi \,dx = 0 \quad \hbox{for any}\  \varphi \in E.
    $$
Thus $g(v_{0})=0$.  
Since $F_{p_0+1} \subset C(\R^{N}\setminus\{0\})$ and $g(s)\not= 0$ in a neighborhood 
of $0$ but excluding $0$ by (g1'), we conclude $v_{0}\equiv 0$ also for $N=2$.

\medskip

By the compactness of the embedding $F_{p_0+1}\subset L_{loc}^r(\R^N)$ for
$r\in [1,2^*)$, Step 1 implies 
    $$  v_j\to 0\quad \hbox{strongly in}\ L_{loc}^r(\R^N)
        \quad \hbox{for}\ r\in [1,2^*).
        \eqno\label[4.22]
    $$

\medskip

\noindent
{\sl Step 2. $v_{j} \to 0$ strongly in $L^{p_0+1}(\R^{N})$ and conclusion
}

\smallskip

\noindent
Since $(v_{j})_{j=1}^{\infty}$ is bounded in $F_{p_0+1}$, by Lemma \ref[Lemma:4.10] 
there exist constants $\delta$, $R_{0}> 0$ such that
    $$  \abs{v_{j}(x)} \leq \delta \quad \hbox{for all}\  \abs{x} \geq R_{0}\
        \hbox{and}\ j\in \N.
    $$
By (g1'), we may assume that for a constant $C'>0$
    $$  g(s)s \geq C' {s}^{p_0+1} \quad \hbox{for}\  s\in [0,\delta].
    $$
Thus
    $$  \eqalign{
        C' \, &\norm{v_{j+}}_{L^{p_0+1}(\abs{x}\geq R_0)}^{p_0+1}
        \leq \int_{\abs{x}\geq R_0} g(v_{j+})v_{j+} \, dx \cr
        &= \int_{\R^N} g(v_j)v_j \, dx - \int_{\R^N} g(-v_{j-})(-v_{j-}) \, dx 
        - \int_{\abs{x}\leq R_0} g(v_{j+})v_{j+} \, dx\cr
        &= (I) - (II) - (III). \cr}
    $$
We observe $(I) \to 0$ by \ref[4.20], 
$(II) \to 0$ by \ref[4.21], 
$(III) \to 0$ by \ref[4.22].
Thus $\norm{v_{j+}}_{L^{p_0+1}(\abs{x}\geq R_0)} \to 0$.
We compute
    $$  \eqalign{
        \norm{v_j}_{p_0+1}^{p_0+1} 
        &= \norm{v_{j+}}_{L^{p_0+1}(\abs{x}\geq R_0)}^{p_0+1}
        + \norm{v_{j+}}_{L^{p_0+1}(\abs{x}\leq R_0)}^{p_0+1}
        + \norm{v_{j-}}_{p_0+1}^{p_0+1} \cr
        &= o(1) + (I) + (II). \cr}
    $$
We observe $(I) \to 0$ by \ref[4.22], $(II) \to 0$ by \ref[4.21].  
Therefore $v_j \to 0$ strongly in $L^{p_0+1}(\R^N)$, 
which is a contradiction to \ref[4.18].  
Thus $(u_{j})_{j=1}^{\infty}$ is bounded in $F_{p_0+1}$.  
Thus Lemma \ref[Lemma:4.11] is proved. \QED

\medskip

\noindent
{\bf
c) Convergence of $(u_j)$ and conclusion
}

\medskip

\noindent
By Lemma \ref[Lemma:4.11], we may assume $u_{j}\wlimit u_{0}$ weakly in $F_{p_0+1}$
after extracting a subsequence.

\medskip

\proclaim Lemma \label[Lemma:4.12]. 
$u_{0}=0$ and $u_{j}\to 0$ strongly in $L^{p_0+1}(\R^{N})$.

\medskip

\claim Proof.
Clearly $u_{0}\in F_{p_0+1}$ is a non-negative solution of \ref[4.15].  
Corollary \ref[Corollary:4.9] implies $u_{0}=0$.
By \ref[4.4] and \ref[4.5],
    $$  \int_{\R^{N}} G(u_{j}) - {N-2\over 2N} g(u_{j})u_{j} dx \to 0. 
        \eqno\label[4.23]
    $$
By (g1') we have for some $\delta, C>0$
    $$  G(s) - {N-2\over 2N} g(s)s \geq C s^{p_0+1}
        \quad \hbox{for}\ s \in [0,\delta].
    $$
On the other hand, by the boundedness of $(u_{j})_{j=1}^{\infty}$ 
in $F_{p_0+1}$, by Lemma \ref[Lemma:4.10] there exists $R_{0}>0$ such that
    $$  \abs{u_{j}(x)} \leq \delta \quad 
        \hbox{for all}\  \abs{x}\geq R_{0}\ \hbox{and}\ j\in \N.
    $$
Thus
    $$  \eqalign{
    C&\norm{u_{j+}}_{L^{p_0+1}(\abs{x}\geq R_{0})}^{p_0+1}
    \leq \int_{\abs x\geq R_0} G(u_{j+}) - {N-2\over 2N} g(u_{j+})u_{j+}\, dx \cr
    &= \int_{\R^N} G(u_j) - {N-2\over 2N} g(u_j)u_j\, dx
        - \int_{\R^N} G(u_{j-}) - {N-2\over 2N} g(-u_{j-})(-u_{j-})\, dx \cr
    &\quad - \int_{\abs x\leq R_0} G(u_{j+}) - {N-2\over 2N} g(u_{j+})u_{j+}\, dx \cr
    &= (I)-(II)-(III).\cr}
    $$
We observe $(I)\to 0$ by \ref[4.23], $(II)\to 0$ by \ref[4.6], 
$(III)\to 0$ by the compactness of the embedding $F_{p_0+1}\subset L^r(\abs x\leq R_0)$.
Thus $\norm{u_{j+}}_{L^{p_0+1}(\abs{x}\geq R_{0})}\to 0$.

As in Step 2 of the proof of Lemma \ref[Lemma:4.11], we conclude that $u_j\to 0$ strongly 
in $L^{p_0+1}(\R^N)$.
\QED

\medskip

\claim End of the proof of Lemma \ref[Lemma:4.6].
Since $u_j\wlimit 0$ weakly in $F_{p_0+1}$ and the embedding 
$F_{p_0+1}\subset L^q(\R^N)$ is compact for $q\in (p_0+1,2^*)$,
Lemma \ref[Lemma:4.12] implies that $u_j\to 0$ strongly
in $L^r(\R^N)$ for $r\in [p_0+1,2^*)$.  Thus
    $$  \int_{\R^{N}} g(u_{j})u_{j}\,dx \to 0.
    $$
On the other hand, by \ref[4.5], we have 
    $   \int_{\R^{N}} g(u_{j})u_{j} dx \to Nb
    $,
which is a contradiction when $b\not=0$.
Thus $\mu_{j}\to 0$ cannot occur for a $(PSPC)_b$ sequence 
$(\lambda_{j},u_{j})$ with $b\not=0$. 
\QED

\medskip

\claim Proof of Proposition \ref[Proposition:4.1].  
Proposition \ref[Proposition:4.1] follows from 
Lemmas \ref[Lemma:4.2], \ref[Lemma:4.4], \ref[Lemma:4.5] and
\ref[Lemma:4.6].
\QED

\medskip

\BS{\label[Section:5]. Application to cubic-quintic type equations and $L^2$-supercritical problems}

\BSS{\label[Subsection:5.1]. Existence results}
We study the cubic-quintic problems and $L^2$-supercritical problems 
as applications of our Theorems \ref[Theorem:1.4] and \ref[Theorem:1.6].  
We assume (g0), (g1'), (g2).

We study here nonlinearities which are $L^2$-supercritical at $s=0$, that is,
we assume that (g1') holds with
    $$  1+{4\over N} < p_0 < 2^* - 1.
    $$
For the behavior as $s \to \infty$, we consider two cases:  

\smallskip

\item{(1)} $g(s)$ is $L^2$-supercritical but Sobolev subcritical at 
$s = \infty$ (Theorem \ref[Theorem:5.1]);

\item{(2)} $g(s)$ is $L^2$-subcritical at $s = \infty$,
which includes cubic-quintic type nonlinearities.  
(Theorem \ref[Theorem:5.4], Remark \ref[Remark:5.5]).

\smallskip

\noindent
In this subsection we give the existence results whose special cases are Theorem \ref[Theorem:1.1] and B.
Proofs of our existence results will be given in the following subsections.

We start with (1). We have the following result:  

\medskip

\proclaim Theorem \label[Theorem:5.1].
Assume (g0), (g1'), (g2) and (p1).
Moreover assume
{\parindent=1.5\parindent
\item{(a1)} there exists $\theta \in (2 + {4\over N}, 2^*)$ and $L > 1$ such that
    $$  0 < \theta G(s) \leq g(s) s \quad \hbox{for all}\ s \geq L.
    $$
}
Then for any $m > 0$, $(*)_m$ has at least one solution 
$(\mu,u) \in (0,\mu_*) \times H^1(\R^N)$ with $u(x) > 0$ in $\R^N$.

\medskip

\claim Remark \label[Remark:5.2].
{\sl 
We note ${N\over N-2} \le 1 + {4\over N}$ for $N \ge 4$.
Thus (p1) holds if and only if  one of the following conditions holds

\smallskip

{\parindent=1.8\parindent
\item{(N2)} $N=2$ and $p_0\in (3,\infty)$;
\item{(N3a)} $N=3$ and $p_0\in ({7\over 3}, 3]$;
\item{(N3b)} $N=3$, $p_0\in (3,5)$ and \ref[4.1];
\item{(N4)} $N\geq 4$, $p_0\in (1+{4\over N}, 2^*-1)$ and \ref[4.1].

}

\smallskip

\noindent
We observe that (p0) in Theorems \ref[Theorem:1.1] and \ref[Theorem:1.2] is nothing but (N2) or (N3a).  
Thus our Theorem \ref[Theorem:1.2] is a special case of Theorem \ref[Theorem:5.1].

} 

\medskip

\claim Remark \label[Remark:5.3].  
{\sl 
(a1) holds if $\lim_{s\to\infty} {g(s)\over s^\nu}=1$ for some 
$\nu\in (1+{4\over N},2^*-1)$.
Conversely (a1) implies the existence of a constant $C>0$ such that
    $$  G(s) \geq C s^{\theta}, \quad g(s) \geq C s^{\theta-1}
        \quad \hbox{for}\ s \geq L.     \eqno\label[5.1]
    $$
} 

\negthree

This type of existence result has been well studied since the pioneering work 
of Jeanjean [\cite[J:24]],  
who considered the constrained problem \ref[1.1] and applied the MP minimax method 
on the $L^2$-sphere $\calS_m$.  The results in [\cite[J:24]] has been improved in 
[\cite[BS:5], \cite[BM:8], \cite[JL1:26]],
in which they assume some global conditions to show the existence results.
For example,

\smallskip

\bitem ((f4) in [\cite[JL1:26]])  $\widetilde{G}(s)$ is strictly increasing on $(0,\infty)$, 
where $\widetilde{G}(s) = g(s)s - 2 G(s)$.  

\bitem ((f5) in [\cite[JL1:26]]) When $N \geq 3$,
    $$  g(s)s < 2^*  G(s) \quad \hbox{for all}\  s \in \R \setminus \{0\}.
        \eqno\label[5.2]
    $$

\smallskip

\noindent
We emphasize that Theorem \ref[Theorem:5.1] under (N2) or (N3a) 
(e.g. Theorem \ref[Theorem:1.2]) gives an existence result
without such global conditions.

\medskip

Next we consider situation (2) and assume that $g(s)$ is 
$L^2$-subcritical at $s = \infty$.

\medskip

\proclaim Theorem \label[Theorem:5.4].
Assume (g0), (g1'), (g2) and (p1).
Moreover assume one of the following two conditions:
{
\itemitem{(a2-1)} $\mu_* < \infty$, where $\mu_*$ is defined in \ref[1.4].
When $N=2$, assume moreover \m
$\limsup_{s\to\infty} {g(s)\over s^{1+{4\over N}}} \leq 0$;
\itemitem{(a2-2)} $\mu_* = \infty$ and there exists 
$\theta \in (2, 2 + {4\over N})$ and $L > 1$ such that
    $$  \theta G(s) \geq g(s) s \quad \hbox{for}\ s \geq L.   \eqno\label[5.3]
    $$
}
Then
{
\item{(i)} $c_* \equiv \inf_{\mu \in (0, \mu_*)} c_-(\mu) > 0$;
\item{(ii)} if $m = c_*$, $J_m(\mu,u)$ has at least one critical point $(\mu_0,u_0)$,
that is, $(*)_m$ has at least one positive solution.  Moreover $J_m(\mu_0,u_0)>0$. 
\item{(iii)} if $m \in (c_*, \infty)$, $J_m(\mu,u)$ has at least two critical points 
$(\mu_1, u_1)$ and $(\mu_2, u_2)$,  
that is, $(*)_m$ has at least two distinct positive solutions. Moreover
    $$  J_m(\mu_1,u_1) > 0, \quad J_m(\mu_1,u_1) > J_m(\mu_2,u_2).
    $$
}

\noindent
We have the following remarks on conditions (g2) and (a2-2).

\medskip

\claim Remark \label[Remark:5.5] {(c.f. Berestycki-Lions [\cite[BL:6]])}.
{\sl
Let $s_0>0$ be a number given in (g3). If there exists a number $s_1>s_0$
such that $g(s_1)=0$, we set
    $$  \widetilde g(s) = \cases{
            g(s)    & $s \in (-\infty,s_1]$,\cr
            0       & $s \in (s_1,\infty)$. \cr}
    $$
Since solutions of $-\Delta u + \mu u = \widetilde g(u)$ in $\R^N$
satisfy $u(x)\le s_1$ for all $x \in \R^N$ by the maximal principle,
we may replace $g$ with $\widetilde g$.  
We note that $\mu_*\in (0,\infty)$ for $\widetilde G(s)=\int_0^s g(\tau)\, d\tau$
and (a2-1) holds for the truncated problem. 
Choosing the smallest $s_1>s_0$ with $g(s_1)=0$, we may assume $G(s_1)>0$ and
thus $\widetilde g(s)$ satisfies \ref[4.1] if $g(s)$ satisfies \ref[4.1].

For example, when $N=3$, $g(s)=s^3-s^5$ does not satisfy condition (g2),
however considering the truncated function $\widetilde g(s)$ of $g(s)$,
we can apply Theorem \ref[Theorem:5.4].
} 

\medskip

\claim Remark \label[Remark:5.6].
{\sl 
(a2-2) holds if $\lim_{s\to\infty} {g(s)\over s^\nu}=1$ for some 
$\nu\in (1,1+{4\over N})$.
Conversely, (a2-2) implies
    $$  \eqalignno{
        &G(s) > 0 \quad \hbox{in } [L,\infty), &\label[5.4]\cr
        &\hbox{there exists $C>0$ such that $G(s) \leq C s^{\theta}$, 
            $g(s) \leq C s^{\theta-1}$ for $s \in [L,\infty)$}. &\label[5.5]\cr}
    $$
In fact, since $\mu_* = \infty$, there exists a sequence $s_j \to \infty$ 
such that $G(s_j)>0$ and $G(s_j) \to \infty$.  
Thus we have
    $$  {g(s)\over G(s)} \leq {\theta \over s} \quad \hbox{in a neighborhood of}\ s_j,
    $$
which implies that for $s$ in a neighborhood of $s_j$ with $s < s_j$,
    $$  \log G(s) \geq \log G(s_j) + \theta \log {s\over s_j} 
        \geq \log G(s_j) + \theta \log {L\over s_j}. 
        \eqno\label[5.6]
    $$
We can see that \ref[5.6] holds for $s\in [L,s_j)$ as long as $G(s)>0$ and 
we find $G(s) > 0$ for all $s \in [L, s_j]$.  
Since $s_j \to \infty$, we have \ref[5.4].
\ref[5.5] follows from \ref[5.3] and \ref[5.4].

} 

\medskip

\noindent
Theorem \ref[Theorem:1.1] is a special case of Theorem \ref[Theorem:5.4].  
Here we give a proof of Theorem \ref[Theorem:1.1].

\medskip

\claim Proof of Theorem \ref[Theorem:1.1].
The conditions (g0), (g1'), (p1) clearly hold in the setting of Theorem \ref[Theorem:1.1].
Under condition (2) in Theorem \ref[Theorem:1.1], conditions (g2) and (a2-2) hold 
(see also Remark \ref[Remark:5.6]) and Theorem \ref[Theorem:5.4] is applicable.  
Under condition (1) in Theorem \ref[Theorem:1.1], we have 
$\mu_* <\infty$.  If $g(s)$ satisfies $g(s)\geq 0$ for 
large $s>0$, clearly we have $\lim_{s\to \infty} {g(s)\over s^\theta}=0$ for any $\theta>1$
and thus (g2) and (a2-1) hold and Theorem \ref[Theorem:1.1] is applicable.  
When there exists a sequence $(t_n)_{n=1}^\infty$ such that $t_n\to\infty$ and $g(t_n)<0$,
we argue as in Remark \ref[Remark:5.5] and we can reduce $(*)_m$ to 
 $-\Delta u + \mu u = \widetilde g(u)$ in $\R^N$, where $\widetilde g$ is a trucation
of $g$ and $\widetilde g(s)=0$ holds for large $s>0$.  
Thus Theorem \ref[Theorem:5.4] is applicable after truncation.
\QED

\medskip

Results related to Theorem \ref[Theorem:5.4] are obtained by Jeanjean-Lu [\cite[JL3:28]] and
Theorem \ref[Theorem:5.4] is largely motivated by [\cite[JL3:28]].
They consider a class of nonlinearities $g(s)$ which is $L^2$-subcritical 
at $s=\infty$ and $L^2$-supercritical at $s=0$, which includes 
cubic-quintic type nonlinearities.
Among other results they study the existence of two solutions under a global condition:

\smallskip

\bitem there exists $\theta \in (2, 2^*)$ such that
    $$  \theta G(s) \geq g(s) s \quad \hbox{for all}\ s \in \R.
        \eqno\label[5.7]
    $$

\smallskip

\noindent
They show that there exists $m_{**} > 0$ such that for $m > m_{**}$, 
$(*)_m$ has two distinct solutions.  
Their result can be applied to
    $$  g(s) = a \abs{s}^{p-2}s - b \abs{s}^{q-2}s \eqno\label[5.8]
    $$
where $a,b>0$ and $2 + {4\over N} < p < q < {2N\over N-2}$ ($N \geq 3$),  
$2 + {4\over N}< p < q < \infty$ ($N=1,2$).  
In particular, they can treat the cubic-quintic problem ($N=3, p=4, q=6$).  
In our Theorem \ref[Theorem:1.1], which is a special case of Theorem \ref[Theorem:5.4],
we do not require global conditions like \ref[5.7].

Our solutions are obtained as applications of Theorem \ref[Theorem:1.4] and 
Theorem \ref[Theorem:1.6] and properties of $b_m(\mu)$ are important;
$(\mu_1, u_1)$ (resp. $(\mu_2, u_2)$) 
corresponds to a local maximum (resp. local minimum) of $b_m(\mu)$.  

\bigskip
\bigskip

\line{
\hfil\pdffig{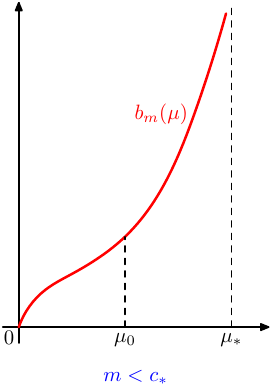}
\hfil\pdffig{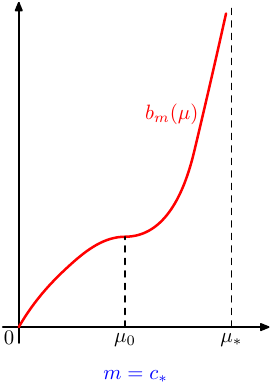}
\hfil\pdffig{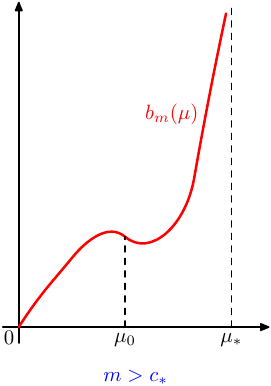}
\hfil}

\medskip

\centerline{Graphs of $b_m(\mu)$}


\bigskip

The approach taken in [\cite[JL3:28]] is quite different from our approach; 
a solution is obtained as a minimizer of the minimizing problem:
    $   \inf \{ \calI(u); \; u \in \calS_m \}
    $
and another solution $(\mu_2, u_2)$ is obtained as a local minimum in
    $   \calS_m^\rho = \{ u \in \calS_m ; \, \norm{ \nabla u }^2_{L^2(\R^N)} > \rho \}
    $
for a suitable $\rho > 0$.

\medskip

\claim Remark \label[Remark:5.7].  
{\sl 
(i) It is natural to ask whether $c_* = m_{**}$ holds or not. 
Here $m_{**}$ is introduced in [\cite[JL3:28]]. \m
(ii) For $m \in (0, c_*)$, we observe
    $(b_m)'_+(\mu) = c_-(\mu) - m  \geq c_* - m > 0$ for all $ \mu\in (0,\mu_*)$.
By Remark \ref[Remark:1.5] (i), any least energy solution $u$ of $-\Delta u+\mu u =g(u)$,
$\mu\in (0,\mu_*)$ is not a solution of $(*)_m$.\m
(iii)
In [\cite[KOPV:32], Section \raw[2]], Killip, Oh, Pacovnicu and Visan study 
the structure of positive radial solutions of $(*)_m$ precisely when
$g(s) = s^3 - s^5$, $N = 3$.  
They find $\mu_* = {3\over 16}$ in our notation and they show that there exists 
a unique positive radial solution $u_\mu$ of  $(**)$ for each $\mu \in (0,\mu_*)$
applying the results of Serrin-Tang ([\cite[ST:37], Theorems \raw[1', 4]]).
Thus $c_-(\mu) = c_+(\mu)=\half\norm{u_\mu}_2^2$ for each $\mu \in (0,\mu_*)$. 
They also establish numerous results for $u_\mu$ including the precise asymptotic 
behavior of $c_-(\mu)$ as $\mu\to 0^+$ and $\mu\to\mustarm$.\m
(iv)
In [\cite[KOPV:32], \raw[Conjecture 2.3]], they conjecture that there exists 
$\mu_0 \in (0,\mu_*)$ such that $\mu \mapsto c_-(\mu)$
is strictly decreasing in $(0,\mu_0)$ and strictly increasing
in $(\mu_0,\mu_*)$.  
Their conjecture is supported by their numerical results.
If their conjecture is valid, then $c_* = c_-(\mu_0)$ and the number of positive 
solutions of $(*)_m$ is
$0$ (resp.\ $1$, $2$) for
$m \in (0,c_*)$ (resp.  $m = c_*$,  $m \in (c_*,\infty)$).

} 

\medskip


\BSS{\label[Subsection:5.2]. Proof of Theorem \ref[Theorem:5.1]}
Since (a1) implies \ref[5.1], recalling $\theta\in(2+{4\over N},2^*)$,  
we have $\mu_* = \infty$ and assumption (ii) of Proposition \ref[Proposition:3.1] holds.
We also note that assumption (ii) of Corollary \ref[Corollary:3.4] follows 
from (g1') with $p_0\in (1+{4\over N}, 2^*-1)$.
Thus we have
    $$  \eqalignno{
        &a(\mu) \to 0, \quad {a(\mu)\over \mu} \to \infty
            \quad \hbox{as}\ \mu \to 0^+,                    &\label[5.9]\cr
        &{a(\mu)\over \mu} \to 0 \quad \hbox{as } \mu \to \infty. &\label[5.10] \cr}
    $$
We also note that (b3) in Proposition \ref[Proposition:4.1] is nothing but (a1).
Thus $\widetilde J_m$ satisfies $(PSPC)_b$ for all $b \neq 0$.

\medskip

\claim Proof of Theorem \ref[Theorem:5.1].
By \ref[5.9]--\ref[5.10],
    $$  b_m(\mu) = a(\mu) - m \mu = \mu \left( {a(\mu)\over \mu} - m \right)
    $$
satisfies
    $$  \eqalign{
        &b_m(\mu) \to 0\quad  \hbox{as}\  \mu \to 0^+, \cr
        &\hbox{there exists $\delta>0$ such that $b_m(\mu) > 0$ in $(0,\delta]$}, \cr
        &b_m(\mu) \to -\infty\quad \hbox{as}\  \mu \to \infty. \cr}
    $$
Thus, $b_m(\mu):(0,\infty) \to \R$ has a positive global maximum at some 
$\mu_0 \in (0,\infty)$.
Since $\widetilde J_m$ satisfies the $(PSPC)_b$ condition 
for all $b \geq b_m(\mu_0)>0$,
Theorem \ref[Theorem:1.6] implies that $J_m$ has a critical point 
$(\mu_\#, u_\#)\in (0,\infty)\times E$ such that
    $$  J_m(\mu_\#, u_\#) \geq b_m(\mu_0) > 0,
    $$
that is, a positive solution of $(*)_m$.
\QED

\medskip

\BSS{\label[Subsection:5.3]. Proof of Theorem \ref[Theorem:5.4]}
We recall Remark \ref[Remark:5.6].  
Since \ref[5.5] follows from (a2-2), assumption (i) of 
Proposition \ref[Proposition:3.1] holds under (a2-2). Thus
    $$  {a(\mu)\over \mu} \to \infty, \qquad c_-(\mu) \to \infty 
            \quad \hbox{as}\ \mu \to \infty. \eqno\label[5.11] 
    $$
Under (a2-1), i.e., $\mu_* < \infty$, by Proposition \ref[Proposition:3.2] we have
    $$  a(\mu) \to \infty, \qquad c_-(\mu) \to \infty 
        \quad \hbox{as}\ \mu \to \mustarm. \eqno\label[5.12]
    $$
Assumption (ii) of Corollary \ref[Corollary:3.4] also holds 
as in Theorem \ref[Theorem:5.1] and thus
    $$  a(\mu) \to 0, \quad {a(\mu)\over \mu} \to \infty, 
        \quad c_-(\mu) \to \infty \quad \hbox{as}\ \mu \to 0^+, \eqno\label[5.13]
    $$
Since (b1) (resp. (b2)) in Proposition \ref[Proposition:4.1] is nothing 
but (a2-1) (resp. (a2-2)), 
$(PSPC)_b$ holds for $\widetilde{J}_m(\lambda,u)$ for all $b\not= 0$.

\medskip

\claim Proof of Theorem \ref[Theorem:5.4]. \m
We deal with two cases $\mu_* = \infty$ and $\mu_* < \infty$ 
simultaneously.  In what follows,  
$\mu \to \mu_*$ means $\mu \to \infty$ (when $\mu_* = \infty$) or 
$\mu \to \mustarm$ (when $\mu_* < \infty$).

(i) By \ref[5.11]--\ref[5.12], we have $c_-(\mu) \to \infty$ as $\mu \to 0^+$ and 
$\mu \to \mu_*$.
Since $\mu \mapsto c_-(\mu)$ is lower semi-continuous and $c_-(\mu) > 0$ for all 
$\mu \in (0,\mu_*)$, we have
    $c_* \equiv \inf_{\mu \in (0,\mu_*)} c_-(\mu) > 0$,
and it is attained at some $\mu_0 \in (0,\mu_*)$.

(ii) Suppose $m = c_* (= c_-(\mu_0))$. 
Let $u_0 \in \calC(\mu_0)$ with $\half\norm{u_0}_2^2=c_-(\mu_0)=c_*$.  
Then $u_0$ is a least energy solution of
$-\Delta u_0 + \mu_0 u_0 = g(u_0)$ with $\half\norm{u_0}_2^2 = c_*$.  
Clearly $(\mu_0,u_0)$ is a solution of $(*)_m$ with $m = c_*$.
We also note that $(b_m)_\pm'(\mu) = a_\pm'(\mu)-m = c_\mp(\mu)-m\geq c_*-m=0$
for all $\mu\in (0,\mu_*)$.  Thus $b_m(\mu)$ is non-decreasing in $(0,\mu_*)$.
In particular, $J_m(\mu_0,u_0)=b_m(\mu_0)>0$.

(iii) 
By \ref[5.11]--\ref[5.13],
    $   b_m(\mu) = a(\mu) - m\mu = \mu \left( {a(\mu)\over \mu} - m \right)
    $
satisfies
    $$  \eqalignno{
        &b_m(\mu) \to 0\quad  \hbox{as}\  \mu \to 0^+, &\label[5.14]\cr
        &\hbox{there exists $\epsilon > 0$ such that $b_m(\mu) > 0$ in $(0,\epsilon]$}, 
            &\label[5.15]\cr
        &b_m(\mu) \to \infty \quad \hbox{as}\  \mu \to \mu_*. &\label[5.16]\cr}
    $$
Now suppose $m > c_* = c_-(\mu_0)$. Then by Theorem \ref[Theorem:1.3] we have
    $$  (b_m)'_+(\mu_0) = a'_+(\mu_0) - m = c_-(\mu_0) - m = c_* - m < 0.
    $$
Thus there exists small $\delta>0$ with $\mu_0+\delta<\mu_*$ such that
    $$  b_m(\mu_0+\delta) < b_m(\mu_0). 
    $$
Consider two intervals $[0, \mu_0 + \delta]$, $[\mu_0, \mu_*)$ and set
    $$  \overline{b} = \max_{\mu \in (0, \mu_0 + \delta]} b_m(\mu), \quad
        \underline{b} = \min_{\mu \in [\mu_0, \mu_*)} b_m(\mu).
    $$
We have
    $$  0<\overline b, \quad \underline b<\overline b.
    $$
and we can easily see that $\overline b$ is a topological local maximum of $J_m$
in $(0,\mu_0+\delta)$ and $\underline b$ is a local minimum of $J_m$ in 
$(\mu_0,\mu_*)$. 

Thus by Theorem \ref[Theorem:1.4], $\underline{b}$ is a critical value of $J_m(\mu,u)$.
Since $(PSPC)_b$ holds for $b \geq \overline{b}>0$, we can apply Theorem \ref[Theorem:1.6] to show
the existence of a critical value $b_\# \geq \overline b$.
Since $\underline{b} < b_\#$, $J_m(\mu,u)$ has at least two positive solutions
corresponding to $\underline b$ and $b_\#$.
\QED

\BS{\label[Section:6]. Proofs of abstract results}
In this section we give proofs to Theorems \ref[Theorem:1.3], \ref[Theorem:1.4] 
and \ref[Theorem:1.6].

\BSScap{\label[Subsection:6.1]. Proof of Theorem \ref[Theorem:1.3] when $N\geq 3$}{\ref[Subsection:6.1]. Proof of Theorem \ref[Theorem:1.3] when N geq 3}
To show Theorem \ref[Theorem:1.3], we use ideas from [\cite[JT1:29]].  
Let $\mu\in (0,\mu_*)$.  In [\cite[JT1:29]], an optimal path for the MP value 
$a(\mu)$ is given explicitly.  When $N \geq 3$, the optimal path is of the form
    $$  \gamma(\tau) = u(x/(T\tau)) \in \Gamma_\mu,
    $$
where $u\in \calC(\mu)$ and $T\gg 1$.


For an arbitrary given compact interval $[\mu_1,\mu_2] \subset (0,\mu_*)$, 
we will show the right differentiability of $a(\mu)$ at $\mu_0 \in (\mu_1,\mu_2)$.
As in Section \ref[Subsection:2.2], we use the notation
    $$  g_\mu(s) = g(s) - \mu s, \quad 
        G_\mu(s) = G(s) - {\mu\over 2}s^2   \eqno\label[6.1]
    $$
and we set
    $$  \widetilde{\calC} = \{ (\mu,u) \,;\, 
        \mu \in [\mu_1,\mu_2], \, u \in \calC(\mu) \}.      \eqno\label[6.2]
    $$
Note that $\widetilde{\calC}$ is compact.

First we claim that there exist $h_0>0$ and $T \gg 1$ such that
    $$  I(\mu+h,u(x/T)) < 0 \quad \hbox{for all}\  
        (\mu,u)\in \widetilde{\calC}, \, \abs{h}\le h_0. \eqno \label[6.3]
    $$
In fact, by the Pohozaev identity ${N-2\over 2}\norm{\nabla u}_2^2 -N\intRN G_\mu(u)\, dx=0$, 
we have for $(\mu,u)\in \widetilde{\calC}$ and $h\in\R$,
    $$  \eqalignno{
        I(\mu+h,u(x/T))
        &= {T^{N-2}\over 2}\norm{\nabla u}_2^2
            + T^N \left({h\over 2}\norm{u}_2^2 - \int_{\R^N} G_\mu(u)\,dx \right), \cr
        &= {T^{N-2}\over 2}\norm{\nabla u}_2^2
            + T^N\left({h\over 2}\norm{u}_2^2 
            - {N-2\over 2N}\norm{\nabla u}_2^2\right). &\label[6.4] \cr}
    $$
Since there exist $c_0,c_1>0$ such that
    $$  c_0 \le \norm{u}_2^2, \, \norm{\nabla u}_2^2 \le c_1 
        \quad \hbox{for}\  (\mu,u)\in \widetilde{\calC},
    $$
we find that \ref[6.3] holds for small $h_0>0$ and large $T\gg 1$.
In particular, $u(x/(T\tau)) \in \Gamma_{\mu+h}$ for all 
$(\mu,u)\in \widetilde{\calC}$ and $\abs{h}\le h_0$.

To prove Theorem \ref[Theorem:1.3] for $N\geq 3$, we note that
    $$  \eqalign{
        a(\mu) &= \max_{\tau\in[0,1]} I(\mu,u(x/(T\tau)))
            = \max_{\tau\in[0,T]} I(\mu,u(x/\tau)), \cr
        a(\mu+h) &\le \max_{\tau\in[0,T]} I(\mu+h,u(x/\tau)). \cr}
    $$
For each $h \in [0,h_0]$, we choose $u_h \in \calC(\mu_0+h)$.
For $h=0$, we also assume that $u_0 \in \calC(\mu_0)$ satisfies
    $   \half \norm{u_0}_2^2 = c_-(\mu_0)
    $.
We have
    $$  D_1(h)  \leq a(\mu_0+h) - a(\mu_0) \leq D_2(h),
    $$
where
    $$  \eqalign{
        D_1(h) &= I(\mu_0+h, u_h(x)) - \max_{\tau \in [0,T]} I(\mu_0, u_h(x/\tau)), \cr
        D_2(h) &= \max_{\tau \in [0,T]} I(\mu_0+h, u_0(x/\tau)) - I(\mu_0, u_0(x)).\cr}
    $$
We compute $D_1(h)$. We set
    $$  A_h = {N-2\over 2}\norm{\nabla u_h}_2^2, \quad 
        B_h = N \int_{\R^N} G_{\mu_0+h}(u_h)\,dx, \quad
        C_h = {N\over 2}\norm{u_h}_2^2.
    $$
We have $B_h=A_h$ by Pohozaev identity.  Thus we have
    $$  \eqalignno{
        I(\mu_0+h, u_h(x)) &= {1\over N-2}A_h - {1\over N}B_h 
            = {2\over N(N-2)} A_h\cr
        I(\mu_0, u_h(x/\tau)) &= {1\over N-2} A_h \tau^{N-2}
            - {1\over N}(B_h + hC_h)\tau^N \cr
        &= {1\over N-2} A_h \tau^{N-2}- {1\over N}(A_h + hC_h)\tau^N. \cr}
    $$
We can see that $\tau \mapsto I(\mu_0, u_h(x/\tau))$ takes maximum at 
$\tau=\tau_h= \left(1+{hC_h\over A_h}\right)^{-1/2}$.
Thus
    $$  \eqalign{
    D_1(h) &= {2\over N(N-2)}A_h 
            -\left({1\over N-2} A_h \tau_h^{N-2} - {1\over N}(A_h+hC_h)\tau_h^N\right)\cr
    &= {2\over N(N-2)}\left(A_h-A_h\left(1+{hC_h\over A_h}\right)^{-{N-2\over 2}}\right). \cr}
    $$

Now take a sequence $h_j \to 0^+$ such that
${1\over h_j} D_1(h_j) \to \liminf_{h\to 0^+} {1\over h} D_1(h)$.
By the compactness of $\widetilde\calC$, we may assume that $u_{h_j}(x) \to u_*(x)$ strongly in $E$
for some $u_*(x) \in \calC(\mu_0)$, which implies 
    $$  A_{h_j} \to A_* \equiv {N-2\over 2}\norm{\nabla u_*}_2^2 > 0, \quad
        C_{h_j} \to C_* \equiv {N\over 2}\norm{u_*}_2^2 > 0. 
    $$
Thus
    $$  \eqalign{
        D_1(h_j) &= {2\over N(N-2)}\left(A_{h_j}
        - A_{h_j}\left(1 - {N-2\over 2}{h_j C_{h_j}\over A_{h_j}} + o(h_j)\right)\right) \cr
        &= {1\over N}h_j C_{h_j}  +o(h_j), \cr}
    $$
which implies
    $$  {1\over h_j}D_1(h_j) \to {1\over N} C_* = \half \norm{u_*}_2^2 \geq c_-(\mu_0).
    $$
Thus
    $$  \liminf_{h \to 0^+} {a(\mu_0+h)-a(\mu_0)\over h}
        \geq \liminf_{h\to 0^+} {1\over h}D_1(h) 
        =\lim_{j\to\infty} {1\over h_j} D_1(h_j) 
        \geq c_-(\mu_0). 
    $$
In a similar way, we can show
    $$  \limsup_{h \to 0^+} {a(\mu_0+h)-a(\mu_0)\over h}
        \leq \limsup_{h \to 0^+} {1\over h} D_2(h)
        = \half \norm{u_0}_2^2 = c_-(\mu_0).
    $$
Thus we have
    $$  \lim_{h\to 0^+} {a(\mu_0+h)-a(\mu_0)\over h} = c_-(\mu_0).
    $$
That is, $a(\mu)$ is right-differentiable and
$a'_+(\mu) = c_-(\mu)$.  
In a similar way, we can also show $a'_-(\mu) = c_+(\mu)$.  
\QED

\medskip


\BSS{\label[Subsection:6.2]. Proof of Theorem \ref[Theorem:1.3] when $N=2$}
In [\cite[JT1:29]], an optimal path for the MP value $a(\mu)$ is given explicitly 
also for $N=2$.  
Let $\mu \in (0,\mu_*)$ and $u \in \calC(\mu)$.  
The explicit path in [\cite[JT1:29]] is given by joining the following 3 curves in $E$:
    $$  \eqalign{
        &[0,1] \to E;\ s \mapsto s\,u(x/t_0), \cr
        &[t_0,t_1] \to E;\ t \mapsto u(x/t), \cr
        &[1,s_1] \to E;\ s \mapsto s\,u(x/t_1), \cr}
    $$
where $0 < t_0 < 1 < t_1,\ s_1 > 0$ are suitably chosen constants.  
We will use a related path to prove Theorem \ref[Theorem:1.3].

We will take an arbitrary given compact interval $[\mu_1,\mu_2] \subset (0,\mu_*)$ and
we will show the right differentiability of $a(\mu)$ at
$\mu_0 \in (\mu_1,\mu_2)$.

As in Section \ref[Subsection:2.1] we set $\widetilde\calC$ by \ref[6.2] and
use notation \ref[6.1].
For any $(\mu,u) \in \widetilde{\calC}$, we have
    $$  \int_{\R^2} G_\mu(u)\,dx = 0, \quad
        \norm{\nabla u}_2^2 = \int_{\R^2} g_\mu(u)\,u\,dx.
        \eqno\label[6.5]
    $$
Since $\widetilde{\calC}$ is compact, there exist constants $c_0,c_1>0$ such that
    $$  c_0 \le \norm{\nabla u}_2^2,\ \int_{\R^2} g_\mu(u)\,u\,dx \le c_2
        \quad \hbox{for all}\  (\mu,u)\in \widetilde{\calC}. 
        \eqno\label[6.6]
    $$
We also remark that for $\delta_0>0$ small,
    $$  \eqalignno{
    &\int_{\R^2} G_\mu(su)\,dx \ge \half c_0\abs{s-1} 
        \quad \hbox{for}\  s\in[1,1+\delta_0], &\label[6.7]\cr
    &\int_{\R^2} G_\mu(su)\,dx \le -\half c_0\abs{s-1} 
        \quad \hbox{for}\  s\in[1-\delta_0,1]. &\label[6.8]\cr}
    $$
In fact, by \ref[6.6], there exists $\delta_0>0$ such that 
for $(\mu,u)\in\widetilde{\calC}$,
    $$  {d\over ds}\int_{\R^2} G_\mu(su)\, dx 
        = \int_{\R^2} g_\mu(su)u\,dx \ge \half c_0 \quad 
        \hbox{for}\  s\in[1-\delta_0,1+\delta_0].
    $$
Together with $\int_{\R^2}G_\mu(u)\,d x=0$, we have \ref[6.7]--\ref[6.9].

To prove Theorem \ref[Theorem:1.3], for a given $\epsilon \in (0,1)$, we will find 
an optimal path for $a(\mu)$, $\mu \in (\mu_1,\mu_2)$, joining the following paths:
    $$  \eqalign{
        &L_1 : [0,1-\delta_\epsilon] \to E;\ s \mapsto s u(x/t_\epsilon), \cr
        &L_2 : [t_\epsilon,1-\epsilon] \to E;\ t \mapsto (1-\delta_\epsilon)u(x/t), \cr
        &L_3 : [1-\delta_\epsilon,1] \to E;\ s \mapsto s u(x/(1-\epsilon)), \cr
        &L_4 : [1-\epsilon,1+\epsilon] \to E;\ t \mapsto u(x/t), \cr
        &L_5 : [1,1+\delta_\epsilon] \to E;\ s \mapsto s u(x/(1+\epsilon)), \cr
        &L_6 : [1+\epsilon,\ell_\epsilon] \to E;\ t \mapsto (1+\delta_\epsilon) u(x/t). \cr}
    $$
Here $\delta_\epsilon \in (0,1)$, $t_\epsilon \in (0,1-\epsilon)$, 
$\ell_\epsilon\in (1+\epsilon,\infty)$ are constants 
which depend on $\epsilon$ but not on $(\mu,u)\in \widetilde{\calC}$.
We choose $\delta_\epsilon$, $t_\epsilon$, $\ell_\epsilon$
in the following lemma.

\medskip

\proclaim Lemma \label[Lemma:6.1].  
For any $\epsilon \in (0,1)$ there exist $\delta_\epsilon \in (0,\delta_0)$, 
$t_\epsilon \in (0,1-\delta_\epsilon)$, $\ell_\epsilon\in (1+\epsilon,\infty)$ and
$h_\epsilon>0$ such that for 
$(\mu,u)\in \widetilde{\calC}$ and $\abs{h}\leq h_\epsilon$:
    $$  \eqalignno{
        &{d\over ds} I(\mu+h, s u(x/(1+\epsilon))) < 0
            \quad \hbox{for}\ s \in [1,1+\delta_\epsilon],
            &\label[6.9] \cr
        &{d\over ds} I(\mu+h, s u(x/(1-\epsilon))) > 0. 
            \quad \hbox{for}\ s \in [1-\delta_\epsilon,1],
            &\label[6.10] \cr
        &{d\over dt} I(\mu+h, (1+\delta_\epsilon) u(x/t)) < 0 
            \quad \hbox{for}\  t \in [1+\epsilon,\infty), &\label[6.11] \cr
        &{d\over dt} I(\mu+h, (1-\delta_\epsilon) u(x/t)) > 0 
            \quad \hbox{for}\  t \in (0,1-\epsilon], &\label[6.12]\cr
        &I(\mu+h, (1+\delta_\epsilon) u(x/\ell_\epsilon)) <0, 
            &\label[6.13]\cr
        &{d\over ds} I(\mu+h, s u(x/t_\epsilon)) > 0 
            \quad \hbox{for}\  s \in (0,1).     &\label[6.14]\cr}
    $$

\medskip

\claim Proof.  
First we show \ref[6.9]. We can show \ref[6.10] in a similar way.
    $$  \eqalignno{
    &{d\over ds} I(\mu+h, s u(x/(1+ \epsilon))) \cr
    &= {d\over ds} \left\{ {s^2\over 2}\norm{\nabla u}_2^2 
        + (1+\epsilon)^2\left({h\over 2} s^2 \norm{u}_2^2 
        - \int_{\R^2} G_\mu(su)\,dx\right) \right\} \cr
    &= s\norm{\nabla u}_2^2 + (1+ \epsilon)^2\left( h s \norm{u}_2^2 
        - \int_{\R^2} g_\mu(su) u\,dx \right). &\label[6.15] \cr}
    $$
As $h \to 0,\ s \to 1$, we have by \ref[6.5] uniformly 
in $(\mu,u)\in\widetilde{\calC}$,
    $$  
    \ref[6.15] \to \norm{\nabla u}_2^2 - (1+\epsilon)^2 \int_{\R^2} g_\mu(u)u\,dx
    = (1-(1+\epsilon)^2)\norm{\nabla u}_2^2 <0.
    $$
Thus there exist $\delta_\epsilon\in(0,1)$, $h'_\epsilon>0$ such that  
\ref[6.9] holds for $s\in[1,1+\delta_\epsilon]$ and $\abs{h}\leq h'_\epsilon$.
In a similar way we can show \ref[6.10].

Next we show \ref[6.11].  We compute
    $$  {d\over dt} I(\mu+h,(1+\delta_\epsilon)u(x/t)) 
        = 2t\left({(1+\delta_\epsilon)^2\over 2} h\norm{u}_2^2 
            - \int_{\R^2} G_\mu((1+\delta_\epsilon)u)\,dx\right). \eqno\label[6.16]
    $$
By \ref[6.7], we have 
    $   \int_{\R^2} G_\mu((1+\delta_\epsilon)u)\,dx \ge \half c_0\delta_\epsilon
    $.
Thus there exists $h_\epsilon\in(0,h'_\epsilon)$ independent of 
$(\mu,u)\in\widetilde{\calC}$ such that 
    $$  {(1+\delta_\epsilon)^2\over 2}h\norm{u}_2^2 
        - \int_{\R^2} G_\mu((1+\delta_\epsilon)u)\,dx 
        \leq -{1\over 4}c_0\delta_\epsilon<0
        \quad \hbox{for}\  \abs{h}\le h_\epsilon.
        \eqno\label[6.17]
    $$
Thus we have \ref[6.11].  
We can show \ref[6.12] in a similar way. 
It follows from \ref[6.16] and \ref[6.17] that
    $$  I(\mu+h, (1+\delta_\epsilon) u(x/t)) \to -\infty
        \quad \hbox{as}\ t\to\infty\ \hbox{uniformly in}\ (\mu,u)\in\widetilde\calC.
    $$
Thus there exists $\ell_\epsilon\in (1+\epsilon,\infty)$, for which \ref[6.13] holds.

Finally we show \ref[6.14].  We compute
    $$
    {d\over ds} I(\mu+h, s u(x/t))
    = s\left(\norm{\nabla u}_2^2 + t^2\Big(h\norm{u}_2^2 
        - \int_{\R^2}{g_\mu(su)\over su}u^2\,dx\Big)\right). 
    $$
Since $\norm{u}_\infty$ is uniformly bounded for $(\mu,u)\in\widetilde{\calC}$ and
${g_\mu(\tau)\over \tau}\to -\mu$ as $\tau\to 0$, 
there exists $C>0$ such that
    $$  \pnorm{g_\mu(su)\over su}_\infty \le C \quad 
        \hbox{for all}\ (\mu,u)\in\widetilde{\calC}\ \hbox{and}\ s\in(0,1].
    $$
Thus,
    $$  \pabs{h\norm{u}_2^2 - \int_{\R^2}{g_\mu(su)\over su}u^2\,dx}
        \le (\abs{h}+C)\norm{u}_2^2 \leq C'.
    $$
Therefore for small $t=t_\epsilon \in (0,1-\epsilon)$ we have \ref[6.14].  
\QED

\medskip


For $\epsilon \in (0,1)$, we choose and fix $\delta_\epsilon \in (0,\delta_0)$, 
$t_\epsilon \in (0,1-\delta_\epsilon]$, 
$\ell_\epsilon \in (1+\epsilon,\infty)$ and 
$h_\epsilon > 0$ by Lemma \ref[Lemma:6.1].  
For $(\mu,u) \in \widetilde{\calC}$, we consider a path joining 6 paths $L_1$--$L_6$
and we denote it by $\gamma(u;\tau):[0,1]\to E$ after reparametrization.
By Lemma \ref[Lemma:6.1] we observe that  
for $\abs{h} \leq h_\epsilon$, 
$I(\mu+h, \gamma(u;\tau))$ 
is strictly increasing along $L_1$--$L_3$ and 
strictly decreasing along $L_5$--$L_6$.  
Along $L_4$, we have by \ref[6.5]
    $$  {d \over d\tau} I(\mu+h, u(x/\tau)) 
        = {d \over d\tau} \Bigl( \half \norm{\nabla u}_2^2 
            + {\tau^2 h \over 2}\norm{u}_2^2 \Bigr)
        = \tau h \norm{u}_2^2
        \ \cases{
            >0 & for $h>0$, \cr
            =0 & for $h=0$, \cr
            <0 & for $h<0$. \cr} 
    $$
Thus $I(\mu+h, \gamma(u;\tau))$ takes its maximum  
at the end point of $L_4$ if $h>0$,  
at the starting point of $L_4$ if $h<0$.  
We also remark that $\gamma(u;\tau)$ is an optimal path for $a(\mu)$
if $h=0$.

Thus, for $h \in (0,h_\epsilon)$
    $$  \eqalignno{
    &\max_{\tau \in [0,1]} I(\mu+h, \gamma(u;\tau))
    = I(\mu+h, u(x/(1+\epsilon))) \cr
    &= \half \norm{\nabla u}_2^2 
        + (1+\epsilon)^2 \left( {h\over 2}\norm{u}_2^2 
        - \int_{\R^2} G_\mu(u) \,dx\right) 
    = \half \norm{\nabla u}_2^2 
        + {h\over 2} (1+\epsilon)^2
\norm{u}_2^2 \cr
    &= a(\mu) + {h\over 2}(1+\epsilon)^2\norm{u}_2^2. &\label[6.18]\cr}
    $$
Similarly,
    $$  \max_{\tau \in [0,1]} I(\mu-h, \gamma(u;\tau))
        = I(\mu-h, u(x/(1-\epsilon))) 
        = a(\mu) - {h\over 2}(1-\epsilon)^2 \norm{u}_2^2. \eqno\label[6.19]
    $$

\medskip

\claim End of the proof of Theorem \ref[Theorem:1.3] for $N=2$.  \m
Let $\mu_0 \in (\mu_1,\mu_2)$. For each $h \in [0,h_\epsilon]$ we choose  
$u_h \in \calC(\mu_0+h)$.  
For $h=0$, we also assume $u_0 \in \calC(\mu_0)$ satisfies  
$\half \norm{u_0}_2^2 = c_-(\mu_0)$.  

Then we have  
    $$  D_1(h) \leq a(\mu_0+h) - a(\mu_0) \leq D_2(h),
    $$
where
    $$  \eqalign{
        D_1(h) &= I(\mu_0+h, u_h(x)) 
            -\max_{\tau \in [0,1]} I(\mu_0, \gamma(u_h;\tau)), \cr
        D_2(h) &= \max_{\tau \in [0,1]} I(\mu_0+h, \gamma(u_0;\tau)) 
            - I(\mu_0, u_0(\tau)). \cr}
    $$
Applying \ref[6.19] to $\mu=\mu_0+h$,
    $$  D_1(h) = a(\mu_0+h) 
            - \left(a(\mu_0+h) - {h\over 2}(1-\epsilon)^2 \norm{u_h}_2^2\right)
        =  {h\over 2}(1-\epsilon)^2 \norm{u_h}_2^2.
    $$
By \ref[6.18],
    $$  D_2(h) = \left(a(\mu_0) + {h\over 2}(1+\epsilon)^2 \norm{u_0}_2^2\right)
                - a(\mu_0)
        = {h\over 2}(1+\epsilon)^2 \norm{u_0}_2^2.
    $$
We assume $\liminf_{h\to 0^+}\norm{u_h}_2^2=\lim_{j\to\infty} \norm{u_{h_j}}_2^2$
for some subsequence $h_j\to 0^+$.  Moreover
we may assume that $u_{h_j} \to u_*$ strongly in $E$ for some $u_*\in\calC(\mu_0)$.
Thus we have
    $$  \eqalign{
    (1-&\epsilon)^2 c_-(\mu_0) 
    \leq {(1-\epsilon)^2 \over 2} \norm{u_*}_2^2 
    =\liminf_{h\to 0^+} {D_1(h)\over h} \cr
    &\leq \liminf_{h\to 0^+} {a(\mu_0+h)-a(\mu_0)\over h}
    \leq \limsup_{h\to 0^+} {a(\mu_0+h)-a(\mu_0)\over h} \cr
    &\leq \limsup_{h\to 0^+} {D_2(h)\over h}
    = {(1+\epsilon)^2 \over 2}\norm{u_0}_2^2
    = (1+\epsilon)^2 c_-(\mu_0). \cr}
    $$
Since $\epsilon>0$ is arbitrary, we have  
$a'_+(\mu_0) = c_-(\mu_0)$. 
Similarly we can show  
$a'_-(\mu_0) = c_+(\mu_0)$.
\QED


\medskip


\BSS{\label[Subsection:6.3]. Proof of Theorem \ref[Theorem:1.4]}
Suppose that $b_m(\mu)$ takes a local minimum at $\mu_0 \in (0,\mu_*)$.  
Then we have
    $$  (b_m)'_+(\mu_0) \geq 0, \quad (b_m)'_-(\mu_0) \leq 0. \eqno\label[6.20]
    $$
Recall that $b_m(\mu) = a(\mu) - m\mu$ and thus
    $   (b_m)'_\pm(\mu) = a'_\pm(\mu) - m = c_\mp(\mu) - m
    $.
Since $c_-(\mu) \leq c_+(\mu)$, \ref[6.20] implies
    $$  c_+(\mu_0) = c_-(\mu_0) = m. \eqno\label[6.21]
    $$
Thus, we have $(b_m)'_+(\mu_0) = (b_m)'_-(\mu_0) = 0$ and 
$b_m(\mu)$ is differentiable at $\mu_0$.

\ref[6.21] implies
    $$  \half\norm{u}_2^2 = m \quad \hbox{for all}\  u \in \calC(\mu_0).
    $$
Thus, any $u \in \calC(\mu_0)$ solves
    $  -\Delta u + \mu_0 u = g(u)$, $\half\norm{u}_2^2 = m$,
that is, $(\mu_0, u)$ is a solution of $(*)_m$.  
\QED

\medskip

\BSS{\label[Subsection:6.4]. Proof of Theorem \ref[Theorem:1.6]}
We use an idea from [\cite[CGIT:15]] and we will observe that 
a topologically non-trivial local maximum of $b_m(\mu)$ 
induces a local linking structure of $J_m(\mu,u)$ in $(0,\infty)\times E$.

First we set up 2-dimensional minimax method for $J_m$.

\medskip

\noindent
{\bf 
a) 2-dimensional minimax method for $J_m(\mu,u)$
}

\medskip

\noindent
We take an interval $[\mu_1,\mu_2]\subset (0,\infty)$ and 
a point $\mu_0\in (\mu_1,\mu_2)$ such that
    $$  b_m(\mu_0) > \max\{ b_m(\mu_1), b_m(\mu_2)\}. \eqno\label[6.22]
    $$
We set
    $$  D = [\mu_1,\mu_2]\times [0,1].
    $$
In Section {b}), we will define 
$\chi_0(\mu,\tau)\in C(D,(0,\infty)\times E)$
and we will consider the following minimax value:
    $$  \widehat{b} = \inf_{\sigma\in\Lambda} \max_{(\mu,\tau)\in D} 
        J_m(\sigma(\mu,\tau)), \eqno \label[6.23]
    $$
where
    $$  \Lambda = \{\sigma(\mu,\tau)\in C(D,(0,\infty)\times E); 
        \ \sigma(\tau,\mu)=\chi_0(\tau,\mu)\ 
        \hbox{for}\  (\mu,\tau)\in \partial D\}. \eqno \label[6.24]
    $$
We will show that $\widehat{b}$ is a critical value of $J_m$ in Section {c}).

\medskip

\noindent
{\bf
b) Choice of $\chi_0(\mu,\tau)$
}

\medskip

\noindent
We choose $\chi_0\in C(D,(0,\infty)\times E)$ with the following properties:
    $$  \eqalignno{
    &\chi_0(\mu,0) = (\mu,0) \quad \hbox{for}\  \mu\in [\mu_1,\mu_2], &\label[6.25]\cr
    &\chi_0(\mu,1) = (\mu,e_0) \quad \hbox{for}\  \mu\in [\mu_1,\mu_2], &\label[6.26]\cr
    &\chi_0(\mu_i,\tau) = (\mu_i,\zeta_i(\tau)) \quad 
        \hbox{for}\  \tau\in [0,1], \ i=1,2, &\label[6.27]\cr}
    $$
where $e_0\in E$, $\zeta_i(\tau)\in C([0,1],E)$ satisfy
    $$  \eqalignno{
    &I(\mu,e_0)<0 \quad \hbox{for}\  \mu\in [\mu_1,\mu_2], &\label[6.28]\cr
    &I(\mu_i,\zeta_i(\tau)) \leq a(\mu_i) \quad 
        \hbox{for}\  \tau\in [0,1],\ i=1,2, &\label[6.29]\cr
    &\zeta_i(0)=0,\quad \zeta_i(1)=e_0 \quad \hbox{for}\  i=1,2. &\label[6.30]\cr}
    $$
First we take an optimal path $\zeta_2(\tau)\in \Gamma_{\mu_2}$ 
for the MP value $a(\mu_2)$.  Such a path is given in [\cite[JT1:29]] explicitly
(see also Sections \ref[Subsection:6.1] and \ref[Subsection:6.2]).  
Setting $e_0=\zeta_2(1)$, we have
    $$  I(\mu,e_0) \leq I(\mu_2,e_0)<0 \quad \hbox{for all}\  \mu\in [\mu_1,\mu_2].
    $$
Thus $e_0\in \{u;\, I(\mu_1,u)<0\}$. 
Recalling that $\{u; I(\mu_1,u)<0\}$ is path-connected, there exists an optimal path
$\zeta_1(\tau)\in \Gamma_{\mu_1}$ for $a(\mu_1)$ with $\zeta_1(1)=e_0$.
We observe that $e_0$, $\zeta_1$, $\zeta_2$ satisfy \ref[6.28]--\ref[6.30].

Setting $\chi_0(\mu,\tau)$ for $(\mu,\tau)\in \partial D$ 
by \ref[6.25]--\ref[6.27] and extending continuously onto $D$, we define
    $   \chi_0(\mu,\tau): D \to (0,\infty)\times E
    $.
We also have
    $$  \eqalign{
        &J_m(\chi_0(\mu,0)) = I(\mu,0)-m\mu = -m\mu
            \quad \hbox{for}\ \mu\in [\mu_1,\mu_2], \cr
        &J_m(\chi_0(\mu,1)) = I(\mu,e_0)-m\mu < -m\mu 
            \quad \hbox{for}\ \mu\in [\mu_1,\mu_2], \cr
        &J_m(\chi_0(\mu_i,\tau)) = I(\mu_i,\zeta_i(\tau))-m\mu_i 
            \leq a(\mu_i)-m\mu_i = b_m(\mu_i)\ \hbox{for}\ \tau\in[0,1]\ \hbox{and}\ i=1,2.\cr}
    $$
Moreover under the assumption \ref[6.22] we have
    $$  \max_{ (\mu,\tau)\in \partial D} J_m(\chi_0(\mu,\tau))
        < b_m(\mu_0).           \eqno\label[6.31]
    $$

\medskip


\noindent
{\bf c) Linking property of $\Lambda$
}

\medskip

\noindent
$\Lambda$ enjoys the following linking property.

\medskip

\proclaim Lemma \label[Lemma:6.2].  
For any $\sigma\in \Lambda$, 
    $$  \sigma(D)\cap (\{\mu_0\}\times \calP(\mu_0))\neq \emptyset,
    $$
where $\calP(\widehat \mu)$ is the Pohozaev manifold defined in \ref[2.5].

\medskip


\claim Proof.
The proof consists of several steps.
First we give a property of the Pohozaev functional $P(\mu,u)$ around $u=0$.

\smallskip

\noindent
{\sl
Step 1:
For any $\mu_0 > 0$ there exists $\delta> 0$ such that
    $$  P(\mu_0,u) > 0 \quad \hbox{for}\  0 < \norm{u}_E < \delta.
    $$
}

\vskip -5mm

\noindent
Since
$P(\mu_0,u) = {N-2\over 2} \norm{\nabla u}_2^2 + {N\over 2}\mu_0 \norm{u}_2^2
+ o(\norm{u}^2_E)$ as  $u \sim 0$,
Step 1 holds clearly for $N \geq 3$. For $N=2$, we rely on Proposition \ref[Proposition:2.2].  
By Proposition \ref[Proposition:2.2], we have
    $   \delta \equiv \inf_{u \in \calP(\mu_0)} \norm{\nabla u}_2^2 > 0
    $.
Thus, $u\not\in \calP(\mu_0)$ for $0<\norm u_E<\delta$.  That is,
    $$  P(\mu_0,u) \neq 0 \quad \hbox{for}\  0 < \norm{u}_E < \delta.
    $$
It is easy to find $u$ with small $\norm u_E$ satisfying $P(\mu_0, u)>0$.
Since $\{u; 0 < \norm{u}_E < \delta\}$ is connected, Step 1 holds.

\smallskip

Next for a given $\sigma = (\sigma_1,\sigma_2) \in \Lambda$, 
we give a linking result.

\smallskip

\noindent
{\sl
Step 2: For $\epsilon > 0$ small, there exists 
$z_\epsilon = (\mu_\epsilon,\tau_\epsilon) \in D$ such that
    $$  \sigma_1(z_\epsilon) = \mu_0, 
        \quad P(\mu_0, \sigma_2(z_\epsilon)) = -\epsilon.
    $$
}

\vskip -5mm

\noindent
We use Brouwer degree theory to show Step 2. We define
$F(\mu,\tau): D \to \R^2$ by 
    $$  F(\mu,\tau) = (F_1(\mu,\tau), F_2(\mu,\tau)) 
        = (\sigma_1(\mu,\tau), P(\mu_0,\sigma_2(\mu,\tau))).
    $$
We observe
    $$  \eqalign{
        &F_1(\mu_1,\tau)=\mu_1, \quad F_1(\mu_2,\tau)=\mu_2 
            \quad \hbox{for}\ \tau\in [0,1], \cr
        &F_2(\mu,0)=P(\mu_0,0)=0, \cr
        &F_2(\mu,1)=P(\mu_0,e_0)=NI(\mu,e_0)-\norm{\nabla e_0}_2^2 < 0
             \quad \hbox{for}\ \mu\in [\mu_1,\mu_2].\cr}
    $$
Choosing $\epsilon\in (0,-\max_{\mu\in [\mu_1,\mu_2]} F_2(\mu,1))$, 
we can see that
    $   \deg(F,D,(\mu_0,-\epsilon))\neq 0
    $.
Thus Step 2 follows.

\smallskip

\noindent
{\sl
Step 3: Conclusion
}

\smallskip

\noindent
Let $z_\epsilon \in D$ be a point obtained in Step 2. By Step 1, we have
    $   \norm{\sigma_2(z_\epsilon)}_E \geq \delta
    $.
By the compactness of $D$, we may assume $z_\epsilon \to z_0 \in D$ 
as $\epsilon \to 0$.  $z_0 \in D$ satisfies
    $$  \sigma_1(z_0) = \mu_0, \quad P(\mu_0,\sigma_2(z_0)) = 0, \quad
        \norm{\sigma_2(z_0)}_E \geq \delta.
    $$
Thus $\sigma(z_0) \in \{\mu_0\} \times\calP(\mu_0)$ and 
Lemma \ref[Lemma:6.2] is proved. \QED

\medskip


As a corollary to Lemma \ref[Lemma:6.2] we have

\medskip

\proclaim Proposition \label[Proposition:6.3].
    $$  \widehat{b} \geq b_m(\mu_0).
    $$

\claim Proof. 
By Lemma \ref[Lemma:6.2] and Proposition \ref[Proposition:2.2], we have for any $\sigma\in \Lambda$,
    $$  \eqalignno{
    \max_{(\mu,\tau)\in D} J_m(\sigma(\mu,\tau)) 
    &\geq \inf_{u\in \calP(\mu_0)} J_m(\mu_0,u)
        = \inf_{u\in \calP(\mu_0)} I(\mu_0,u)-m\mu_0 \cr
    &= a(\mu_0)-m\mu_0 = b_m(\mu_0).
        &\label[6.32] \cr}
    $$
Since $\sigma\in \Lambda$ is arbitrary, we have
    $\widehat{b} \geq b_m(\mu_0)$. 
\QED

\medskip


\noindent
{\bf 
d) Proof of Theorem \ref[Theorem:1.6]}

\medskip

\claim End of the proof of Theorem \ref[Theorem:1.6].  
Under \ref[6.22], we show that $\widehat{b}$ defined in \ref[6.23]--\ref[6.24] 
is a critical value of $J_m(\mu,u)$.

Introducing a change of variable $\mu=e^\lambda$, we define
    $$  \eqalign{
        &\widetilde{D} = [\lambda_1,\lambda_2]\times [0,1], \
            \quad \mu_i=e^{\lambda_i}, \ i=1,2, \cr
        &\widetilde{\Lambda} 
            = \{\ \widetilde{\sigma}(\lambda,\tau)\in C(\widetilde{D},\R\times E);\, 
            \widetilde{\sigma}(\lambda,\tau) = \widetilde{\sigma}_0(\lambda,\tau)\ 
            \hbox{for}\  (\lambda,\tau)\in \partial \widetilde{D}\ \}. \cr}
    $$
Here we write $\chi_0(\mu,\tau) = (\chi_{01}(\mu,\tau),\chi_{02}(\mu,\tau))$ and
    $\widetilde{\sigma}_0(\lambda,\tau) 
        = (\log \chi_{01}(e^\lambda,\tau),\ \chi_{02}(e^\lambda,\tau))
    $.

We observe
    $$  \widehat{b} = \inf_{\widetilde{\sigma}\in \widetilde{\Lambda}} 
            \max_{(\lambda,\tau)\in \widetilde{D}} 
            \widetilde{J}_m(\widetilde{\sigma}(\lambda,\tau)).
    $$
By Proposition \ref[Proposition:6.3] and \ref[6.31], we also have
    $$  \widehat{b} \geq b_m(\mu_0) > \max_{(\lambda,\tau)\in \partial\widetilde{D}} 
            \widetilde{J}_m(\widetilde{\sigma}_0(\lambda,\tau)).
    $$
Since $\widetilde{J}_m$ satisfies the $(PSPC)_b$ condition for 
$b\geq b_m(\mu_0)$, in particular, $\widetilde J_m$ satisfies 
the $(PSPC)_{\widehat{b}}$ condition.
Thus we can apply Proposition \ref[Proposition:2.9] to show that $\widehat{b}$ is a critical value of 
$\widetilde{J}_m$, that is, a critical value of $J_m$. 
\QED


\bigskip

\noindent
{\bf Acknowledgments}\m
S.~C. is partially supported by INdAM-GNAMPA.
S.~C. thanks acknowledges financial support from PNRR MUR project PE0000023 NQSTI 
- National Quantum Science and Technology Institute (CUP H93C22000670006). \m
M. G. is supported by 
INdAM-GNAMPA Projects
\lq\lq{}Metodi variazionali per problemi dipendenti da operatori frazionari isotropi e anisotropi'', 
codice CUP E5324001950001, and 
\lq\lq{}Aspetti geometrici e qualitativi di equazioni ellitiche e paraboliche'', codice CUP E53C25002010001. \m
K.~T. is partially supported by JSPS KAKENHI Grant Numbers 
JP22K03380, 24H00024 and 25K07073.

\bigskip

\bibliography

\medskip

{\parindent=1.5\parindent

\bibitem[AS:1] S. N. Armstrong, B. Sirakov, 
Nonexistence of positive supersolutions of elliptic equations via the maximum principle, 
Comm. Partial Differential Equations 36 (2011), no. 11, 2011--2047.

\bibitem[B:2]
L. Baldelli, 
On the cubic-quintic Schr\"odinger equation, 
Nonlinear Diff. Equ. Appl. 32 (2025), article number 105.

\bibitem[BdV:3] T. Bartsch, S. de Valeriola,
Normalized solutions of nonlinear Schr\"odinger equations, 
Arch. Math. (Basel) 100 (2013), no. 1, 75--83.

\bibitem[BMRV:4]
T. Bartsch, R. Molle, M. Rizzi, G. Verzini, 
Normalized solutions of mass supercritical Schr\"odinger equations with potential, 
Comm. Partial Differential Equations 46 (2021), no. 9, 1729--1756.

\bibitem[BS:5]
T. Bartsch, N. Soave, 
A natural constraint approach to normalized solutions of nonlinear Schr\"odinger 
equations and systems, 
J. Funct. Anal. 272 (2017), no. 12, 4998--5037,
Correction, J. Funct. Anal. 275 (2018), no. 2, 516--521.

\bibitem[BL:6] H. Berestycki, P.-L. Lions, 
Nonlinear scalar field equations. I. Existence of a ground state, 
Arch. Rational Mech. Anal. 82 (1983), no. 4, 313--345.

\bibitem[BDS:7] B. Bieganowski, P. d'Avenia, J. Schino, 
Existence and dynamics of normalized solutions to Schr\"odinger equations 
with generic double-behaviour nonlinearities, 
J. Differential Equations 441 (2025), Paper No. 113489, 36.

\bibitem[BM:8]
B. Bieganowski, J. Mederski, 
Normalized ground states of the nonlinear Schr\"odinger equation with at least mass critical growth, 
J. Funct. Anal. 280 (2021), no. 11, Paper No. 108989, 26.

\bibitem[By:9] J. Byeon,
Singularly perturbed nonlinear Dirichlet problems with a general nonlinearity, 
Trans. Amer. Math. Soc. 362 (2010), no. 4, 1981--2001.

\bibitem[CKS:10] R. Carles, C. Klein, C. Sparber, 
On ground state (in-)stability in multi-dimensional cubic-quintic
Schr\"oodinger equations, 
ESAIM Math. Model. Numer. Anal., 57 (2023), no. 2, 423--443.

\bibitem[CS:11]
R. Carles, C. Sparber, 
Orbital stability vs. scattering in the cubic-quintic Schr\"odinger equation, 
Rev. Math. Phys. 33 (2021), no. 3, Paper No. 2150004, 27.

\bibitem[CL:12]
T. Cazenave, P.-L. Lions, 
Orbital stability of standing waves for some nonlinear Schr\"odinger equations, 
Comm. Math. Phys. 85 (1982), no. 4, 549--561.

\bibitem[CGT1:13] S. Cingolani, M. Gallo, K. Tanaka, 
Multiple solutions for the nonlinear Choquard equation with even or odd nonlinearities, 
Calc. Var. Partial Differential Equations 61 (2022), no. 2, Paper No. 68, 34.

\bibitem[CGT3:14] S. Cingolani, M. Gallo, K. Tanaka, 
Infinitely many free or prescribed mass solutions for fractional Hartree equations and Pohozaev identities, 
Adv. Nonlinear Stud. 24 (2024), no. 2, 303--334.

\bibitem[CGIT:15] S. Cingolani, M. Gallo, N. Ikoma, K. Tanaka, 
Normalized solutions for nonlinear Schr\"odinger equations with $L^2$-critical nonlinearity, 
arXiv:{2410.23733}, to appear in Analysis \& PDE.

\bibitem[CT1:16] S. Cingolani, K. Tanaka, 
Ground state solutions for the nonlinear Choquard equation with prescribed mass, 
In: Ferone, V., Kawakami, T., Salani, P., Takahashi, F. (eds) Geometric Properties for Parabolic and Elliptic PDE's,
23--41, 
Springer INdAM Series, vol 47. Springer, Cham. (2021).

\bibitem[CT2:17] S. Cingolani, K. Tanaka, 
Deformation argument under PSP condition and applications, 
Anal. Theory Appl., Vol. 37, No. 2 (2021), 191--208.

\bibitem[CGM:18] S. Coleman, V. Glaser, A. Martin, 
Action minima among solutions to a class of Euclidean field equations, 
Comm. Math. Phys. 58 (1978), no. 2, 211--221.

\bibitem[DST1:19] S. Dovetta, E. Serra, P. Tilli, 
Uniqueness and non-uniqueness of prescribed mass {NLS} ground states on metric graphs, 
Adv. Math. 374 (2020), 107352, 41.

\bibitem[DST2:20] S. Dovetta, E. Serra, P. Tilli, 
Action versus energy ground states in nonlinear Schr\"odinger equations, 
Math. Ann. 385 (2023), no. 3-4, 1545--1576.

\bibitem[HIT:21] J. Hirata, N. Ikoma, K. Tanaka, 
Nonlinear scalar field equations in {$\R^N$}: mountain pass and symmetric mountain pass approaches, 
Topol. Methods Nonlinear Anal. 35 (2010), no. 2, 253--276.

\bibitem[HT:22] J. Hirata, K. Tanaka, 
Nonlinear scalar field equations with $L^2$ constraint: 
Mountain pass and symmetric mountain pass approaches, 
Adv. Nonlinear Stud., Vol. 9, Issue 2, 263--290 (2019).

\bibitem[IM:23]
N. Ikoma, Y. Miyamoto, 
Stable standing waves of nonlinear Schr\"odinger equations 
with potentials and general nonlinearities, 
Calc. Var. Partial Differential Equations 59 (2020), no. 2, Paper No. 48.

\bibitem[J:24] L. Jeanjean, 
Existence of solutions with prescribed norm for semilinear elliptic equations, 
Nonlinear Anal. 28 (1997), no. 10, 1633--1659.

\bibitem[JLe:25]  L. Jeanjean, T. T. Le, 
Multiple normalized solutions for a Sobolev critical Schr\"odinger equation, 
Math. Ann. 384 (2022), no. 1-2, 101--134.

\bibitem[JL1:26] L. Jeanjean, S.-S. Lu, 
A mass supercritical problem revisited, 
Calc. Var. Partial Differential Equations 59 (2020), no. 5, 174.

\bibitem[JL2:27] L. Jeanjean, S.-S. Lu, 
On global minimizers for a mass constrained problem, 
Calc. Var. Partial Differential Equations 61 (2022), no. 6, Paper No. 214, 18.

\bibitem[JL3:28] L. Jeanjean, S.-S. Lu, 
Normalized solutions with positive energies for a coercive problem and application 
to the cubic-quintic nonlinear Schr\"odinger equation, 
Math. Models Methods Appl. Sci. 32 (2022), no. 8, 1557-1588.

\bibitem[JT1:29] L. Jeanjean, K. Tanaka, 
A remark on least energy solutions in {${\bf R}^N$}, 
Proc. Amer. Math. Soc. 131 (2003), no. 8, 2399--2408.

\bibitem[JT2:30] L. Jeanjean, K. Tanaka, 
A note on a mountain pass characterization of least energy solutions, 
Advanced Nonlinear Studies 3 (2003), 461-471.

\bibitem[JZZ:31]
L. Jeanjean, J. Zhang, X. Zhong, 
A global branch approach to normalized solutions for the Schr\"odinger equation, 
J. Math. Pures Appl. (9) 183 (2024), 44--75.

\bibitem[KOPV:32]
R. Killip, T. Oh, O. Pocovnicu, M. Visan, 
Solitons and scattering for the cubic-quintic nonlinear Schr\"odinger equation on 
$\R^3$, 
Arch. Ration. Mech. Anal. 225 (2017), no. 1, 469--548.

\bibitem[MS:33] J. Mederski, J. Schino, 
Least energy solutions to a cooperative system of Schr\"odinger equations 
with prescribed {$L^2$}-bounds: at least {$L^2$}-critical growth, 
Calc. Var. Partial Differential Equations 61 (2022), no. 1, Paper No. 10, 31.

\bibitem[QZ:34] S. Qi, W. Zou, 
Mass threshold of the limit behavior of normalized solutions to Schr\"odinger 
equations with combined nonlinearities, 
J. Differential Equations 375 (2023), 172--205.

\bibitem[Ru:35] H.-J. Ruppen, 
The existence of infinitely many bifurcating branches, 
Proc. Roy. Soc. Edinburgh Sect. A 101 (1985), no. 3-4, 307--320.

\bibitem[Sc:36] J. Schino, 
Normalized ground states to a cooperative system of Schr\"odinger equations 
with generic {$L^2$}-subcritical or {$L^2$}-critical nonlinearity, 
Adv. Differential Equations 27 (2022), no. 7-8, 467--496.

\bibitem[ST:37] J. Serrin, M. Tang, 
Uniqueness of ground states for quasilinear elliptic equations, 
Indiana Univ. Math. J. 49 (2000), no. 3, 897--923.

\bibitem[Sh:38] M. Shibata, 
Stable standing waves of nonlinear Schr\"odinger equations with a general nonlinear term, 
Manuscripta Math. 143 (2014), no. 1-2, 221--237.

\bibitem[So:39] N. Soave, 
Normalized ground states for the {NLS} equation with combined nonlinearities, 
J. Differential Equations 269 (2020), no. 9, 6941--6987.

\bibitem[St:40] C. A. Stuart, 
Bifurcation from the continuous spectrum in the {$L^2$}-theory of elliptic equations on {${\bf R}^n$}, 
Recent methods in nonlinear analysis and applications (Naples, 1980), 231-300, Liguori, Naples (1981).

\bibitem[TVZ:41] T. Tao, M. Visan, X. Zhang, 
The nonlinear Schr\"odinger equation with combined power-type nonlinearities, 
Comm. Partial Differential Equations 32 (2007), no. 7-9, 1281--1343.

\bibitem[WW:42] J. Wei, Y. Wu, 
Normalized solutions for Schr\"odinger equations with critical Sobolev exponent 
and mixed nonlinearities, 
J. Funct. Anal. 283 (2022), no. 6, Paper No. 109574, 46.

} 

\bye